\documentclass[11pt]{article}

\usepackage{fullwidth}
\usepackage{color}
\usepackage[color,matrix,arrow]{xy}
\usepackage[top=1.5in, bottom=1.5in, left=1in, right=1in]{geometry}
\usepackage{amsgen}
\usepackage{amsmath}
\usepackage{amstext}
\usepackage{amsbsy}
\usepackage{amsopn}
\usepackage{amsfonts}
\usepackage{amssymb}
\usepackage{eepic}
\usepackage{graphicx}
\usepackage{epsf}
\usepackage{pstricks}
\xyoption{all}

\def\Box{\square}
\def\edge{\relbar\joinrel\relbar}

\def\mapright#1{\smash{\mathop{\longrightarrow}\limits^{#1}}}

\def\tra#1{\smash{\mathop{\mid\kern
-1pt\joinrel\relbar\joinrel\relbar}\limits^{*}_{#1}}}
\def\longtra#1{\smash{\mathop{\mid\kern
-1pt\joinrel\relbar\joinrel\relbar\joinrel\relbar}\limits^{*}_{#1}}}
\def\vlongtra#1{\smash{\mathop{\mid\kern
-1pt\joinrel\relbar\joinrel\relbar\joinrel\relbar\joinrel\relbar}\limits^{*}_{#1}}}
\def\vvlongtra#1{\smash{\mathop{\mid\kern
-1pt\joinrel\relbar\joinrel\relbar\joinrel\relbar\joinrel\relbar\joinrel\relbar}\limits^{*}_{#1}}}
\def\vvvlongtra#1{\smash{\mathop{\mid\kern
-1pt\joinrel\relbar\joinrel\relbar\joinrel\relbar\joinrel\relbar\joinrel\relbar\joinrel\relbar}\limits^{*}_{#1}}}
\def\etra#1{\smash{\mathop{\mid\kern
-1pt\joinrel\relbar\joinrel\relbar}\limits_{#1}}}

\def\vlongrightarrow{\relbar\joinrel\longrightarrow}

\def\longmapright#1{\smash{\mathop{\vlongrightarrow}\limits^{#1}}}

\def\Rw{\Rightarrow}
\def\oo{\overline}
\def\un{\underline}

\def\wh{\widehat}

\def\C{{\cal{C}}}

\def\M{{\cal{M}}}

\def\mat{\mbox{Mat}}

\def\SPC{\mbox{SPC}}
\def\ker{\mbox{Ker}\,}

\def\max{\mbox{max}}
\def\min{\mbox{min}}
\def\mindeg{\mbox{mindeg}}

\def\LR{\mbox{LR}}

\def\flatx{\mbox{Fl}}

\def\H{{\cal{H}}}
\def\BoR{{\cal{BR}}}

\def\p{\varphi}

\def\inv{^{-1}}

\def\bi{\begin{itemize}}
\def\ei{\end{itemize}}
\def\beq{\begin{equation}}
\def\eeq{\end{equation}}

\newtheorem{T}{Theorem}[section]
\newcommand{\bt}{\begin{T}}
\newcommand{\et}{\end{T}}
\newcommand{\ftd}{$\square$\end{T}}

\newtheorem{Proposition}[T]{Proposition}
\newcommand{\bp}{\begin{Proposition}}
\newcommand{\ep}{\end{Proposition}}
\newcommand{\fpd}{$\square$\end{Proposition}}

\newtheorem{Lemma}[T]{Lemma}
\newcommand{\bl}{\begin{Lemma}}
\newcommand{\el}{\end{Lemma}}
\newcommand{\fld}{$\square$\end{Lemma}}

\newtheorem{Corol}[T]{Corollary}
\newcommand{\bc}{\begin{Corol}}
\newcommand{\ec}{\end{Corol}}
\newcommand{\fcd}{$\square$\end{Corol}}

\newtheorem{Result}[T]{Result}
\newcommand{\br}{\begin{Result}}
\newcommand{\er}{\end{Result}}
\newcommand{\frd}{$\square$\end{Result}}

\newtheorem{Example}[T]{Example}
\newcommand{\be}{\begin{Example}}
\newcommand{\ee}{\end{Example}}

\newtheorem{Problem}[T]{Problem}
\newcommand{\bq}{\begin{Problem}}
\newcommand{\eq}{\end{Problem}}

\newtheorem{Remark}[T]{Remark}
\newcommand{\brm}{\begin{Remark}}
\newcommand{\erm}{\end{Remark}}

\newcommand{\proof}
   {\par\medbreak\noindent{\bf Proof}.\enspace}

\newcommand{\qed}{
$\Box$
\par\bigbreak}

\normalbaselines
\def\abstract#1{\par\bigskip
\begingroup\small
\baselineskip=12truept
\begin{center}ABSTRACT\end{center}
\par\medskip\par\noindent
\null\hfill\hbox{\vbox{\hsize=5truein\noindent#1}}
\hfill\null\par\endgroup\par}

\title{On the Dowling and Rhodes lattices and wreath products}
\author{{\bf Stuart Margolis, John Rhodes and Pedro V. Silva}}

\date{\today}

\begin{document}
\maketitle

\begin{center}\small
2010 Mathematics Subject Classification: 05B35, 05E45, 20M18, 20M30

\bigskip

Keywords: Dowling lattice, Rhodes lattice, wreath product, Brandt groupoid, matroid, lift matroid, Boolean representable simplicial complex
\end{center}

\abstract{Dowling and Rhodes defined different lattices on the set of triples (Subset, Partition, Cross Section) over a fixed finite group G. Although the Rhodes lattice is not a geometric lattice, it defines a matroid in the sense of the theory of Boolean representable simplicial complexes. This turns out to be the direct sum of a complete matroid with a lift matroid of the complete biased graph over G. As is well known, the Dowling lattice defines the frame matroid over a similar biased graph. This gives a new perspective on both matroids and also an application of matroid theory to the theory of finite semigroups.
We also make progress on an important question for these classical matroids: what are the minimal Boolean representations and the minimum degree of a Boolean matrix representation?
}

\section{Introduction}

Let $G$ be a finite group and $X$ a finite set. A partial partition on $X$ is a partition on a subset $I$ of $X$ (including the empty set).  Given $I \subseteq X$, we denote by $F(I,G)$ the collection of all functions
$f:I \rightarrow G$. The group $G$ acts on the left of $F(I,G)$ by $(gf)(x) = gf(x)$ for $f \in F(I,G), g \in G, x \in I$. Then the element $
Gf$ of the quotient set $F(I,G)/G$ is called a {\em cross section} with domain $I$. An
{\em SPC} (Subset, Partition, Cross Section) over $G$ is a triple $(I,\pi,C)$ where $I$ is a subset of $X$, $\pi$ is a partition of $I$ and $C$ is a collection of $k$ cross sections with domains $\pi_1,\ldots,\pi_k$, where $\pi_1,\ldots,\pi_k$ denote the blocks of $\pi$.

We describe next an alternative formalism, which we use most of the time. Given a partition $\pi$ of $I \subseteq X$ and $f,h \in F(I,G)$, we write
$$f \sim_{\pi} h \quad \mbox{ if $f|_{\pi_i} \in G(h|_{\pi_i})$ for each block $\pi_i$ of $\pi$.}$$
Then $\sim_{\pi}$ is an equivalence relation on $F(I,G)$. If we denote by $[f]_{\pi}$ the equivalence class of $f \in F(I,G)$, then
we can view $\SPC(X,G)$ as the set of triples of the form $(I,\pi,[f]_{\pi})$, where $I \subseteq X$, $\pi$ is a partition of $I$ and $f\in F(I,G)$.

The Dowling \cite{Dowling} and Rhodes \cite{pl, KS} orders on $X$ are lattices that are defined on SPCs. Both arise from problems related to wreath products of the form $G \wr S_{X}$, where $S_{X}$ is the symmetric group on $X$ and have many applications. The Dowling order is defined  by refinement and removal of blocks of partial partitions, whereas the Rhodes order is defined by containment of binary relations. It is the purpose of this note to discuss the relationship between these lattices.

The Dowling lattice is a geometric lattice and thus defines a unique simple matroid known as the Dowling geometry or Dowling matroid. On the other hand, the Rhodes lattice is not even atomic, let alone geometric. Nonetheless, in the proper context, the Rhodes lattice does define a matroid. The context is the theory of Boolean representable simplicial complexes (BRSC) \cite{brsc}. One definition of BRSC is by using an appropriate notion of independence over the 2-element Boolean algebra, generalizing the notion of field representable matroids. An equivalent definition is by looking at so called $c$-independent chains \cite{brsc} in an arbitrary finite lattice with a fixed set of join generators. When this definition is applied to a geometric lattice with its set of atoms as generators, one obtains the unique simple matroid on that lattice. When applied to the Rhodes lattice with its set of join irreducible elements as generators, we obtain the direct sum of the complete matroid on $n$ vertices with the lift matroid on the complete $G$-biased graph on $n$-vertices \cite{Zaslavsky2}. In this context, the Dowling matroid is the frame matroid of a similar biased graph. This leads to a new perspective on both matroids.

This surprising connection is one of the main results of this paper. We emphasize that as for the example of the lift matroid of the complete biased graph, the theory of BRSC allows for representing simplicial complexes, including matroids, by lattices smaller than the lattice of flats. We will look for minimal size lattice representations of both the Dowling and Rhodes matroids.

We now give a historical perspective for both the Dowling and Rhodes lattices. The Dowling lattice was first defined for the multiplicative group of a field \cite{Dowling2} and for a general group in \cite{Dowling}, all of this in the early 1970s. The Dowling lattices and matroids play an important role in matroid theory. They form an important class of matroids in the class of frame matroids of biased graphs. See the bibliography \cite{Zaslavsky}.

The Rhodes lattice was defined (unpublished) in 1968 and in print in \cite{lb2, 2J, KS}. The motivation for this definition came from
the complexity theory of finite semigroups, where it is an essential tool in the formulation of the Presentation Lemma. See \cite{qtheor} for details.
See also \cite{lb2, 2J, pl} for specific applications of this concept.

We now outline the paper. In Section 2, we develop a semigroup theoretic viewpoint on the Dowling lattice $Q_n(G)$, by means of principal left ideals of the wreath product $G \wr PT_{n}$, where $G$ is a finite group and $PT_n$ is the monoid of all partial functions on an $n$-element set. This is very important for finite semigroup theory. The Rhodes semilattice $R_n(G)$ and the Rhodes lattice $\wh{R}_n(G)$ are introduced in Section 3. 
In Section 4, we develop a semigroup and groupoid theoretic viewpoint on the Rhodes lattice, using in particular Brandt groupoids. The class of BRSC is introduced in Section 5, and frame and lift matroids in Section 6.
In Section 7, we perform a comparative study of the lattices $Q_n(1)$ and $R_n(1)$, and of the matroids they define.   The matroids defined by $R_n(G)$ and $\wh{R}_n(G)$ are discussed respectively in Sections 8 and 9, establishing the connection to lift matroids. The comparison between the matroids obtained through the Dowling and Rhodes orderings is made in Section 10.
In Section 11 we look at minimal lattice representations for both the Dowling and Rhodes matroids in the sense of the theory of BRSC. We are completely successful in the case of the trivial group. In Section 12 we look at the corresponding problem for representations by Boolean matrices. Both questions are of great importance for the theory. Some open problems are suggested in Section 13.

\section{The Dowling Lattice}

The order on the Dowling lattice is given by refinement and omission of blocks. That is, given two SPCs $(I,\pi,[f]_{\pi})$ and $(J,\tau,[h]_{\tau})$, we define $(I,\pi,[f]_{\pi}) \leq_D
(J,\tau,[h]_{\tau})$ if:
\bi
\item[{\bf 1)}] $J \subseteq I$,
\item[{\bf 2)}] every block of $\tau$ is a union of blocks of $\pi$,
\item[{\bf 3)}] if $\pi_i$ is a block of $\pi$ contained in $J$, then $f|_{\pi_i} \in G(h|_{\pi_i})$.
\ei
The standard notation for this order is $Q_{n}(G)$, where $|X|=n$. It is known that $Q_{n}(G)$ is a geometric lattice and thus defines a unique simple matroid. This also implies that topologically (via its order complex) $Q_{n}(G)$ is a wedge of $n-2$ dimensional spheres \cite{Dowling}. There is a very large literature on Dowling lattices and related concepts \cite{Zaslavsky}.

For example, let $G= \{1\}$ be the trivial group. then $Q_{n}=Q_{n}(1)$ is isomorphic to the lattice of  all partial partitions on an $n$-element set. This lattice is isomorphic to the lattice of all (full) partitions on an $n+1$ element set. The isomorphism is obtained by adding a new element $0$ to $X$ and sending a partial partition $\pi$ on $X$ to the partition on $X \cup \{0\}$ consisting of all the blocks of $\pi$ plus one new block consisting of $0$ and all the elements of $X$ that are not in any block of $\pi$. It is well known that the symmetric group $S_X$ acts on the lattice of all partitions on $X$ and in fact is the automorphism group of this lattice. More generally, the wreath product $G \wr S_{X}$ acts on $Q_{n}(G)$ and this is one of this lattice's most important attributes.

There are a number of equivalent definitions of the Dowling order. The purpose of this section is to give a natural semigroup theoretic interpretation of the Dowling order $Q_{n}(G)$ and the action of $G \wr S_{X}$ that we believe is both new and useful.

There are a number of important posets related to a monoid $M$ and actions of its group of units $G=U(M)$. In particular, let $\mathcal{L}(M)$ be the poset of principal left ideals of $M$ ordered by inclusion. Then $G$ acts on the right of $\mathcal{L}(M)$ by sending $L \in \mathcal{L}(M)$ to $Lg$ for $g \in G$. For example, let $M=M_{n}(F)$ be the monoid of all $n \times n$ matrices over a field $F$. Elementary linear algebra implies that two matrices generate the same principal left ideal if and only if they have the same row space. It follows that $\mathcal{L}(M)$ is the lattice of subspaces of an $n$-dimensional vector space and the action $G=GL_{n}(F)$ is
the usual action of the general linear group on the corresponding projective spaces.

More relevant to this paper is the case of the monoid $PT_{X}$ of all partial functions on a set $X$ acting on the left of $X$. 
If $P$ is a poset, then its opposite $P^{op}$ is the poset with the same set and the opposite order -- that is, $p \leq q$ in $P$ if and only if $ q \leq p$ in $P^{op}$. 

Let $f:X \rightarrow X$ be a
partial function. Recall that $\ker(f)$, the kernel of $f$, is the partial equivalence relation defined on the domain of $f$ by $(x,y) \in \ker(f)$ if and only if
$f(x)=f(y)$. An elementary calculation shows that two partial functions generate the same principal left ideal if and only if they have the same kernel. It follows that $\mathcal{L}(PT_{X})$ can be viewed as the set of partial partitions on $X$. And inclusion in $\mathcal{L}(PT_{X})$ corresponds to the opposite of the Dowling order $Q_{|X|}$. 
Thus, for example, the maximal element in $\mathcal{L}(PT_{X})$ is the (full) partition consisting of all singleton sets and
the minimal element is the empty partition on the empty set (that is, the kernel of the empty function). This is the opposite of the Dowling order.

The next theorem generalizes this last observation to the case of $Q_{n}(G)$. For this it is convenient to view the group $G \wr S_{X}$, where $X$ is an $n$-element set, as the group of $n \times n$ {\em monomial} matrices over $G$. Recall that such a matrix is an $n \times n$ matrix in which every row and column has exactly one element of $G$ and the rest of the elements equal to $0$. More generally, an $n \times n$ {\em column monomial} matrix over $G$ is an $n \times n$ matrix in which each column has
{\bf at most} one element of $G$ and the rest equal to $0$. It is straightforward to see that the collection of all such matrices is isomorphic to the (left) wreath
product $G \wr PT_{n}$ of $G$ with the monoid of all partial functions on an $n$-element set. Furthermore, the group of units of $G \wr PT_{n}$ is $G \wr S_{n}$. The next theorem states that if we order principal left ideals of the monoid $G \wr PT_{n}$ by reverse inclusion we obtain a lattice isomorphic to $Q_{n}(G)$.

\bt
\label{SPCs}
The poset of principal left ideals of the monoid $G \wr PT_{n}$ is a lattice isomorphic to the opposite of the Dowling lattice $Q_{n}(G)$. Furthermore, the usual action of $G \wr S_n$ on $Q_n(G)$ is equivalent to the action of $G \wr S_n$ considered as the group of units of $G \wr PT_{n}$ on its lattice of principal left ideals.
\et

\proof Let $A$ be an $n \times n$ column monomial matrix over $G$. We define an SPC $S(A)$ over $\{1, \ldots n\}$ and $G$. Let $I$  be the set of indices of the non-zero columns of $A$. Let $\pi$ be the partition on $I$ such that $i,j \in I$ are in the same block if and only if the unique non-zero elements in columns $i$ and $j$ are in the same row of $A$. Define a function $f:I \rightarrow G$ by letting $f(i)$ be the unique non-zero element in column $i \in I$. It is easy to check that $S(A)=(I,\pi,[f]_{\pi})$ is an SPC.
Moreover, every SPC over an $n$-element set and group $G$ arises this way. Indeed, let $(I,\pi,[f]_{\pi})$ be an SPC over $\{1,\ldots n\}$ and $G$.
Define the matrix $A$ such that
column $i$ of $A$ is non-zero if and only if $ i \in I$. If the blocks of $\pi$ are $\{\pi_{1}, \ldots, \pi_{k}\}$, the non-zero rows of $A$ are precisely rows 1 to $k$ and such that for each $1 \leq j \leq k$, row $j$ has entry $f(i)$ in column $i$ if $i$ is in $\pi_{j}$ and 0 otherwise. It is clear that $A$ is a column monomial matrix over
$G$ such that $S(A)=(I,\pi,[f]_{\pi})$.

Let $B$ be an $n \times n$ column monomial matrix over $G$. Let $S(BA)=(J,\tau,h)$. We claim that $S(A)\leq_{D} S(BA)$ in the Dowling order. Clearly if a column of $BA$ is
non-zero then the same column in $A$ is non-zero and thus, $J \subseteq I$. If row $i$ of $B$ is the $0$-vector, then so is row $i$ of $BA$. Otherwise, let row $i$ of $B$ have non-zero entries $g_{j_{1}} \ldots g_{j_{k}}$ in columns $j_{1} \ldots j_{k}$. Then row $i$ of $BA$ equals $g_{j_{1}}R_{j_{1}} + \ldots + g_{j_{k}}R_{j_{k}}$ where
$R_{j}$ denotes row $j$ of $A$. Since we are dealing with column monomial matrices, the non-zero positions of row $i$ of $BA$ are the disjoint union of the non-zero
positions of rows $j_{1}, \ldots , j_{k}$ of $A$. It follows that conditions {\bf 2)} and {\bf 3)} in the definition of the Dowling order are fulfilled. Thus if the principal left ideal of a matrix $C$ is contained in that of the matrix $A$, then $S(A) \leq_{D} S(C)$.

Conversely, assume that $S(A) \leq_{D} S(C)$ for $n \times n$ column monomial matrices over $G$. We define a matrix $B$ such that $C=BA$ and this will ensure that the left ideal generated by $C$ is contained in the principal left ideal generated by $A$. Let $r_{i}$, row $i$ of $C$ be non-zero. Since $S(A) \leq_{D} S(C)$, the non-zero
positions of $r_{i}$ are a disjoint union of the non-zero positions of a unique collection of rows $R_{j_{1}}, \ldots, R_{j_{k}}$ of rows of $C$ and that
$r_{i} = g_{j_{1}}R_{j_{1}} + \ldots + g_{j_{k}}R_{j_{k}}$ for some $g_{j_{1}} \ldots g_{j_{k}}$, a tuple of elements of $G$. We define row $i$ of $B$ to be
$g_{j_{1}} \ldots g_{j_{k}}$ with entry $g_{j_{l}}$ in column $j_{l}$ and entry 0 in all other columns. If row $i$ of $C$ is the $0$-vector, then we let row $i$ of B also be the $0$-vector. Direct multiplication shows that $C=BA$ and this completes the proof. The second statement of the theorem is easily verified.
\qed

The results in Theorem \ref{SPCs} have been used, along with appropriate branching rules for the groups $G \wr S_{n}$ to study the representation theory of the monoid $G \wr PT_{n}$ in \cite{Stein1, Stein2}. This has applications to the representation theory of the group $G \wr S_{n}$ via natural semigroup theoretic representations like those discussed in this section. See \cite{Benrepbook} for background in the representation theory of finite monoids. In the language of semigroup theory \cite{CP}, we have proved that two $n \times n$ column monomial matrix over $G$ are in the same $\mathcal{L}$ class, in the sense of Green's relations, if and only if they have the same SPC.

\section{The Rhodes semilattice and the Rhodes lattice}

Let $X$ be  finite nonempty set and let $G$ be a finite group. We consider now a second partial order on $\SPC(X,G)$
based on containment of sets and partitions:
$(I,\pi,[f]_{\pi}) \leq_{R}$ $(J,\tau,[h]_{\tau})$ if:
\bi
\item[{\bf 1)}] 
$I \subseteq J$, 
\item[{\bf 2)}] 
every block of $\pi$ is contained in a (necessarily unique) block of $\tau$,
\item[{\bf 3)}] 
$[h|_I]_{\pi} = [f]_{\pi}$.
\ei
We denote this poset by $R_X(G)$. 

\bl
\label{glb}
$R_X(G)$ is a $\wedge$-semilattice.
\el

\proof 
Let $(I,\pi,[f]_{\pi}), (J,\tau,[h]_{\tau}) \in \SPC(X,G)$. We define 
a relation $R$ on $I \cap J$ by $(x,y) \in R$ if and only if {\bf 1)} $x$ and $y$ are in the same blocks of both $\pi$ and $\tau$ and
{\bf 2)} $f(x)(h(x))\inv = f(y)(h(y))\inv$. It is straightforward that $R$ is a well-defined equivalence relation on $I \cap J$. Let $\rho$ be the partition on $I \cap J$ induced by $R$. If $\rho_i$ is a block of $\rho$, then $f|_{\rho_i} = g(h|_{\rho_i})$ (taking $g = f(x)(h(x))\inv$ for some $x \in \rho_i$) 
and so $[f|_{I \cap J}]_{\rho} = [h|_{I \cap J}]_{\rho}$. Thus $(I \cap J, \rho, [f|_{I \cap J}]_{\rho}) = (I \cap J, \rho, [h|_{I \cap J}]_{\rho})\in \SPC(X,G)$. We claim that 
\beq
\label{glb1}
(I \cap J, \rho, [f|_{I \cap J}]_{\rho}) = ((I,\pi,[f]_{\pi}) \wedge (J,\tau,[h]_{\tau})).
\eeq

It is clear from the definitions that $(I \cap J, \rho, [f|_{I \cap J}]_{\rho})$ is a common lower bound of both $((I,\pi,[f]_{\pi})$ and $(J,\tau,[h]_{\tau})$. 
Suppose now that $(Y,\theta,[t]_{\theta}) \in \SPC(X,G)$ is also a common lower bound of both $((I,\pi,[f]_{\pi})$ and $(J,\tau,[h]_{\tau})$. Then:
\bi
\item
$Y \subseteq I$ and $Y \subseteq J$,
\item
every block $\theta_i$ of $\theta$ is contained in some block $\pi_j$ of $\pi$ and in some block $\tau_k$ of $\tau$,
\item
$[f|_Y]_{\theta} = [t]_{\theta} = [h|_Y]_{\theta}$.
\ei
It follows that $Y \subseteq I \cap J$. If $\theta_i \subseteq \pi_k \cap \tau_k$, take $g_1,g_2 \in G$ such that $f|_{\theta_i} = g_1(t|_{\theta_i})$ and $h|_{\theta_i} = g_2(t|_{\theta_i})$. Then, for all $x,y \in \theta_i$, we get
$$f(x)(h(x))\inv = g_1t(x)(t(x))\inv g_2\inv = g_1t(y)(t(y))\inv g_2\inv = f(y)(h(y))\inv$$
and so $(x,y)\in R$. Hence $\theta_i$ is contained in some block of $\rho$. Since $(f|_{I\cap J})|_Y = f|_Y$ and $[f|_Y]_{\theta} = [t]_{\theta}$, it follows that 
$(Y,\theta,[t]_{\theta}) \leq_R (I \cap J, \rho, [f|_{I \cap J}]_{\rho})$ and so (\ref{glb1}) holds.
\qed

We call $R_X(G)$ the {\em Rhodes semilattice} on $X$ over $G$. Since $R_X(G) \cong R_Y(G)$ when $|X| = |Y|$, we can use the notation $R_{|X|}(G)$ to denote the Rhodes semilattice $R_X(G)$ up to isomorphism. We may also assume that $X = \underline{n} = \{ 1,\ldots,n\}$ for $n = |X|$.

We begin with the case of the trivial group. Then $R_{n}(1)$ can be identified with the set of all pairs $(I,\pi)$ where
$I \subseteq \underline{n}$ and $\pi$ is a partition on $I$. We sometimes write $SP_n$ for this set and the corresponding lattice.

The Rhodes order $\leq_R$ on $SP_n$ is just set inclusion. Thus $(I,\pi) \leq_{R} (J,\tau)$ if and only if $ I \subseteq J$ and every block of $\pi$ is a subset of a (necessarily unique) block of $\tau$. This is the same as demanding that the equivalence relation $R_{\pi}$ corresponding to $\pi$ be a subset of the equivalence relation $R_{\tau}$ corresponding to $\tau$.

We remark that the Rhodes order is different from the Dowling order, even for the trivial group. Indeed, let $n \geq 2$. Then $(\{1, \ldots, n\}, \{\{1\}, \ldots, \{n\}\})$ is the minimal element of the Dowling lattice $Q_{n}(1)$, but is neither the minimal nor maximal element of the Rhodes order $R_{n}(1)$.

Indeed,
$(\{1, \ldots, n-1\}, \{\{1\}, \ldots, \{n-1\}\}) <_{R} (\{1, \ldots, n\}, \{\{1\}, \ldots, \{n\}\})$ and the maximal element of $R_{n}(1)$ is
$(\{1, \ldots, n\}, \{\{1, \ldots, n\}\})$. We will discuss the connection between $R_{n}(1)$ and $Q_{n}(1)$ in more detail later on.

We describe the meet and join in $R_{n}(1)$. 
In view of (\ref{glb1}), we have
$$(I,\pi) \wedge (J,\tau) = (I \cap J, \rho),$$
where the blocks of $\rho$ are obtained by intersecting the blocks of $\pi$ with the blocks of $\tau$. 

Since $SP_n$ is a $\wedge$-semilattice and has $(\un{n}, \{ \un{n} \})$ as top element, then it becomes a lattice with the determined join:
$$ (I,\pi) \vee (J,\tau) = \wedge \{ (K,\theta) \in SP_n \mid (I,\pi) \leq_R (K,\theta) \mbox{ and } (J,\tau) \leq_R (K, \theta) \}.$$
In a more constructive perspective, we can define the equivalence relation $R$ to be the transitive closure of $R_{\pi} \cup R_{\tau}$. Then
$$(I,\pi) \vee (J,\tau) = (I \cup J, \theta),$$
where $\theta$ is the partition on $I \cup J$ defined by $R_{\theta} = R$.

On the other hand, if $|X| = 1$, then $R_1(G)$ has only two elements (determined by the first component)
and is therefore a lattice.

We can give a necessary and sufficient condition for two elements  of $R_n(G)$ admitting a common upper bound (and therefore a join). With this purpose, we introduce the following construction. Given $\alpha_1,\alpha_2 \in \SPC(n,G)$, say $\alpha_i = (I_i,\pi_i,[f_i]_{\pi_i})$ $(i = 1,2)$, we define a graph $\Gamma(\alpha_1,\alpha_2)$ with an edge coloring as follows:
\bi
\item
the vertex set is $I_1 \cup I_2$; 
\item
two distinct vertices are connected by an edge of color $i \in \{ 1,2\}$ if they belong to the same block of $\pi_i$.
\ei

Given an edge $a \edge b$ of color $i$, we may consider an orientation (say $a \mapright{} b$). The label of this directed edge is then $(f_i(a))\inv f_i(b) \in G$. Note that, by replacing $f_i$ by $f'_i \sim_{\pi_i} f_i$, the label remains unchanged. Now, given an oriented path $P$ of the form
\beq
\label{orpath}
a_1 \mapright{} a_2 \mapright{} \cdots \mapright{} a_{m},
\eeq
we define the label of $P$ to be the product of the labels of its edges, following the orientation of the path. 

A {\em necklace} of $\Gamma(\alpha_1,\alpha_2)$ is a cycle alternating edges of both colors (cf. {\em tie your shoes} in \cite{qtheor}). Given a necklace $N$ of the form
\beq
\label{neck}
a_1 \edge a_2 \edge \cdots \edge a_{2m} \edge a_{2m+1} = a_1,
\eeq
we fix an orientation of the cycle and a basepoint, and define the label of $N$ to be the product of the labels of its edges, following the orientation of the cycle. Of course, the label of the necklace depends on both orientation and basepoint, but is unique up to conjugacy and inversion. In particular, the label being 1 does not depend on neither of these factors. 

\bp
\label{existjoin}
The following conditions are equivalent for $\alpha_1,\alpha_2 \in {\rm SPC}(n,G)$:
\bi
\item[(i)] $\alpha_1$ and $\alpha_2$ admit a common upper bound in $R_n(G)$;
\item[(ii)] every necklace of $\Gamma(\alpha_1,\alpha_2)$ has label 1.
\ei
\ep

\proof
(i) $\Rw$ (ii). Write $\alpha_i = (I_i,\pi_i,[f_i]_{\pi_i})$ for $i = 1,2$. Consider a necklace of the form (\ref{neck}). We may assume that $a_1 \edge a_2$ has color 1. 

Let $(I,\pi,[f]_{\pi})$ be a common upper bound of $\alpha_1$ and $\alpha_2$. Let $j \in \{ 1,\ldots, 2m\}$. Then $a_j$ and $a_{j+1}$ belong to the same block $B_j$ of $\pi_{i_j}$, where $i_j = 1$ if $j$ is odd and $i_j = 2$ otherwise. Since $f_{i_j}|_{B_j} \in Gf|_{B_j}$, we get 
$$(f_{i_j}(a_j))\inv f_{i_j}(a_{j+1}) = (f(a_j))\inv f(a_{j+1}).$$
Multiplying the above equalities for $j = 1,\ldots,2m$, we get
$$(f_1(a_1))\inv f_1(a_2)(f_2(a_2))\inv f_2(a_3)(f_1(a_3))\inv f_1(a_4) \ldots (f_2(a_{2m}))\inv f_2(a_1) = 1,$$
hence the necklace has label 1.

(ii) $\Rw$ (i). Let $\alpha = (I,\pi,[f]_{\pi}) \in \SPC(n,G)$ be defined as follows:
\bi
\item
$I = I_1 \cup I_2$;
\item
the blocks of $\pi$ are the connected components of $\Gamma(\alpha_1,\alpha_2)$ (that is, the equivalence relation $R_{\pi}$ is the join of the equivalence relations $R_{\pi_1}$ and $R_{\pi_2}$);
\item
to define the mapping $f:I \to G$, we start by fixing a basepoint $a_0$ in each connected component $C$ of $\Gamma(\alpha_1,\alpha_2)$, and fixing $f(a_0)$ arbitrarily. Then, for every $a \in C$, we consider an oriented path
$a_0 \mapright{} \cdots \mapright{} a$ and 
we define $f(a)$ to be the product of $f(a_0)$ by the label of this path.
\ei

Next we check that $f$ is well defined with respect to $f_1$ and $f_2$. So suppose that we have an alternative path $a_0 \mapright{} \cdots \mapright{} a$. It suffices to show that every closed path has label 1. Straightforward induction reduces the problem to cycles. Since a monochromatic path of the form (\ref{orpath}) has the same label as the edge $a_1 \mapright{} a_m$, and changing the basepoint if needed, we reduce the problem to the necklace case, and the claim now follows from condition (ii). Moreover, our previous remark on the label of edges implies that $[f]_{\pi}$ is well defined with respect to $[f_1]_{\pi_1}$ and $[f_2]_{\pi_2}$. Thus $\alpha$ is well defined.

Now we must show that $\alpha_i \leq_R \alpha$ for $i = 1,2$. Conditions {\bf 1)} and {\bf 2)} are immediate, we focus our attention on condition {\bf 3)}. Let $A$ be a block of $\pi_i$. We must show that $f|_A \in Gf_i|_A$, which amounts to check that $f(a)(f_i(a))\inv$ is constant for every $a \in A$. Let $a,b \in A$ and let $a_0$ be the basepoint of the connected component of $\Gamma(\alpha_1,\alpha_2)$ containing $a$ (and $b$). Considering a path of the form
$a_0 \mapright{} \cdots \mapright{} a \mapright{} b$, we obtain $f(b) = f(a)(f_i(a))\inv f_i(b)$. Thus $f(a)(f_i(a))\inv = f(b)(f_i(b))\inv$ as required.
\qed

Note that, since $(I,\pi) = ((I_1,\pi_1) \vee (I_2,\pi_2))$ in $R_n(1)$, it follows easily that the element $\alpha$ constructed in the proof of the converse implication is indeed the join of $\alpha_1$ and $\alpha_2$.

%
%
%

Now we get the following corollary:

\bc
\label{latt}
The following conditions are equivalent for every finite nonempty set $X$ and every finite group $G$:
\bi
\item[(i)] $R_X(G)$ is a lattice;
\item[(ii)] $|X| = 1$ or $|G| = 1$.
\ei
\ec

\proof
We had already discussed the cases $|X| = 1$ or $|G| = 1$. If $|X|,|G| > 1$, we can obviously construct a necklace with label $\neq 1$, so the claim follows from Proposition \ref{existjoin}.
\qed

If $|X|,|G| > 1$, the simplest way of turning $R_X(G)$ into a lattice $\wh{R}_X(G)$ is by adding a top element $T$, which will be in particular a common upper bound for any pair of distinct elements. Then $\wh{R}_X(G)$ becomes a lattice with the determined join
$$(\alpha \vee \beta) = \wedge\{ \gamma \in \wh{R}_X(G) \mid (\alpha \leq_R \gamma)\mbox{ and } (\beta \leq_R \gamma) \}.$$
Since $R_X(1)$ is already a lattice, we may also use the notation $\wh{R}_X(1) = R_X(1)$.

We note that $\wh{R}_n(G)$ is contractible, that is the order complex \cite{Bjo} of $\wh{R}_n(G)$ is a contractible space. This is because every element of $\wh{R}_n(G)$ is comparable to the element whose set is the whole set, partition is the partition all of whose classes are singletons and the unique CS on this partition.

For $n>1$, this element is neither the maximal nor the minimal element and it is well known that a lattice with an element like this is contractible. On the other hand it is clear that for $n =1$, the space is contractible.

\section{A semigroup theoretic interpretation of the Rhodes lattice}

In this section we give a semigroup theoretic interpretation of the lattice $\wh{R}_{n}(G)$. We start with $G=\{1\}$, the trivial group.

We first use the language of groupoids, where a groupoid is a small category in which all morphisms are isomorphisms. A groupoid is trivial if there is at most one morphism between any two objects. A groupoid is connected if there is at least one morphism between any two objects.

Let $B(1,X)$ be the unique (up to isomorphism) connected trivial groupoid with object set $X$. Thus for all $i,j \in X$, there is exactly one morphism between $i$ and $j$. A subgroupoid $T'$ of a groupoid $T$ is {\em wide} if every object of $T$ is an object of $T'$. In this language, the collection of all wide subgroupoids with object set $X$ is a poset under inclusion of morphism sets with minimal element the groupoid $1_X$ consisting of all identity morphisms and maximal element $B(1,X)$. It is well known that the lattice of all wide subgroupoids of $B(1,X)$ is isomorphic to the partition lattice on $X$. Indeed if $T$ is a wide subgroupoid of $B(1,X)$, then the relation $\pi_T$ on $X$ defined by
$(i,j) \in \pi_T$ if and only if there is a morphism (necessarily unique) from $i$ to $j$ in $T$ is easily seen to be an equivalence relation. Conversely, given an equivalence relation $\pi$ on $X$, define  $T_{\pi}$ to be the wide subgroupoid of $B(1,X)$ with $Mor(i,j) = \{(i,j)\}$ if $(i,j) \in \pi$ and is empty otherwise. It is immediate from the fact that $\pi$ is an equivalence relation that $T_{\pi}$ is indeed a wide subgroupoid of $B(1,X)$. Furthermore, these are inverse operations that preserve order. In the language of inverse semigroup
theory (to be recalled below), this is Theorem 3.2 of \cite{Jonessemimod}

We now wish to generalize this to the lattice of all subgroupoids of $B(1,X)$.  A clear generalization of what we discussed above shows that there is a 1-1 correspondence between subgroupoids of $B(1,X)$ and pairs $(I,\pi)$, where $I$ is a subset of $X$ and $\pi$ is a partition of $I$. Given such a pair, we can define a wide trivial groupoid with object set $I$ and morphisms determined by the partition $\pi$. Conversely, a subgroupoid $T$ of $B(1,X)$ is a wide trivial subgroupoid of $B(1,I)$, where $I$ is the set of objects of $T$. Therefore, there is a unique partition $\pi$ of $I$ corresponding to $T$ and $T$ is coded by the pair $(I,\pi)$. Thus, the collection of of subgroupoids of $B(1,X)$ is in 1-1 correspondence with the collection of partial partitions on $X$. Let $T$ and $T'$ be subgroupoids of $B(1,X)$ and $(I,\pi)$ and $(J,\tau)$ their subset-partition pairs. Then it is clear that $T$ is a subgroupoid of $T'$ if and only if $(I,\pi) \leq_{R} (J,\tau)$, that is the order on subgroupoids of $B(1,X)$ is precisely the Rhodes order $R_{X}(1)$. 

We generalize this to give an interpretation of $R_{X}(G)$ for an arbitrary finite group $G$. In view of Corollary \ref{latt}, we assume that $|X|,|G| > 1$. Let $B(G,X)$ be the groupoid with object set $X$ and for  all $i,j \in X$, the set of morphisms from $i$ to $j$ is $\{(i,g,j) \mid g \in G\}$.
Multiplication is given by $(i,g,j)(j,h,k)=(i,gh,k)$ and the inverse of $(i,g,j)$ is $(j,g^{-1},i)$. It is well known that $B(G,X)$ is, up to isomorphism, the unique connected groupoid with object set $X$ and such that each group $Mor(i,i), i \in X$ is isomorphic to $G$ (it is easy to see by connectedness that if this last fact is true for one object, then it is true for all objects). The next theorem is the main result of this section. It is clear that the intersection of two trivial subgroupoids of a groupoid is itself a trivial groupoid. Easily constructed examples show that the join of two trivial groupoids need not be a trivial groupoid in the lattice of all subgroupoids of a given groupoid. We define the lattice of trivial subgroupoids of a groupoid $G$ to be the collection of all trivial subgroupoids together with a new top element $\mathcal{T}$. We define the meet of two trivial subgroupoids to be their intersection and the join to be their join as subgroupoids if the join is trivial and $\mathcal{T}$ otherwise.

\bt
\label{tsg}
Let $|X|,|G| > 1$. The poset of all trivial subgroupoids of $B(G,X)$ with a new top element adjoined under inclusion is isomorphic to the Rhodes lattice $\wh{R}_{X}(G)$.
\et

\proof To avoid set theoretic arguments, we will assume that $X = \un{n}$. The interested reader can make the changes for arbitrary set $X$.
Let $T$ be a trivial subgroupoid of $B(G,X)$. Let $I$ be the set of objects of $T$, and as previously, the relation $\{(x,y) \mid Mor(x,y) \neq \emptyset\}$ is an equivalence relation on $I$ that defines a partition $\pi$ of $I$. If $I$ is empty, then there is nothing to do. Let $P=\{p_{1}, \ldots, p_{k}\}$ be a block of $\pi$. Then by definition of $\pi$ and triviality, there is a unique morphism between any two elements of $P$. It follows that the function $f:P \rightarrow G$ defined by $f(p_{i})=g_{i}$ if
$(p_{1},g_{i}, p_{i}) \in T, i =1,\ldots, k$ is well defined and that $f(p_{1})=1$ (since the identity is the only isomorphism at $p_{1} \in T$ due to the existence of inverses and closure under composition). We claim that if we choose another base point $p_i$ in place of $p_{1}$, then the function $f_{i}$ defined analogously belongs to $Gf$.

Indeed, $(p_{i},g,p_{j})=(p_{i},(f(p_i))^{-1},p_{1})(p_{1},f(p_j),p_{j})$ by triviality of $T$ and the definition of $f$. Therefore $f_{i}(p_j)=(f(p_i))^{-1}f(p_j)$ and $f_{i} = (f(p_i))\inv f_i \in Gf_i$. By assigning such a function to each block of $\pi$, we have  associated a well defined SPC, $(I,\pi,[f]_{\pi})$ with $T$.

Conversely, let $(J,\tau,[h]_{\tau})$ be an SPC. Pick a block of $\tau$, $\Theta=\{t_{1}, \ldots, t_{m}\}$. Define $T_{\Theta}$ to be the set of morphisms $\{(t_{i},(h(t_i))^{-1}h(t_j),t_{j}) \mid i,j=1,\ldots, m\}$. It is clear that these elements do not depend on the choice of a representative for $[h]_{\tau}$ and that $T_{\Theta}$ defines a connected trivial groupoid with object set $\Theta$. The disjoint union of such $T_{\Theta}$ over the equivalence classes of $\tau$ defines a trivial groupoid associated to $(J,\tau,[h]_{\tau})$. It is straightforward to show that these two operations are inverses of one another and that the containment order on trivial groupoids is isomorphic to the Rhodes order on $R_{X}(G)$. This completes the proof.
\qed

We also get:

\bc
\label{tsgs}
The poset of all trivial subgroupoids of $B(G,X)$ under inclusion is isomorphic to the Rhodes semilattice $R_{X}(G)$.
\ec

We translate this into the theory of inverse semigroups. Let $T$ be a groupoid. Define the {\em consolidation} of $T$ to be the semigroup $T^{0}$ whose elements are the collection of all morphisms of $T$ together with a new element $0$ that will be the zero element of $T^{0}$. The product of two non-zero elements is their product in $T$ if this is defined and $0$ otherwise. It is easy to check associativity. It is well known that the collection of all such semigroups $T^{0}$ is precisely the collection of primitive inverse semigroups with a zero element. See \cite{Lawsonbook} for the relevant definitions. It is easy to see that the set of trivial subgroupoids of $T$ corresponds to the group-free inverse subsemigroups of $T^{0}$. A semigroup is group-free if all its maximal subgroups are trivial groups. Thus we have the following translation of the previous theorem. $B(G,X)^{0}$ is known as the Brandt inverse semigroup over $X$ with structure group $G$. It is the unique (up to isomorphism) completely 0-simple inverse semigroup with an idempotent set of size $|X| +1$ and maximal subgroup isomorphic to $G$ \cite{Lawsonbook}.

\bt

Let $|X|,|G| > 1$. The lattice of group-free inverse subsemigroups of $B(G,X)^{0}$  containing 0 is isomorphic to the Rhodes lattice $\wh{R}_{X}(G)$.

\et

The three authors of the present paper study the BRSC associated to the lattice $R_{n}(1)$ and more generally to the lattice of all subsemigroups of $B(1,X)$ in
\cite{SubsgpBRSC}.

As a final remark, it is well known that the group $G \wr S_{X}$ acts as a group of automorphisms on the groupoid (or corresponding Brandt semigroup) $B(X,G)$ by automorphisms. Indeed, we can consider the morphisms of $B(G,X)$ to be precisely the $|X| \times |X|$ matrices over $G$ that have a non-zero entry in precisely one position by sending the morphism $(i,g,j)$ to the matrix $|X| \times |X|$ that has $g$ in position $(i,j)$ and $0$ elsewhere. Product of matrices when non-zero give composition in $B(G,X)$. Therefore $B(G,X)$ can be considered to be a subset of the monoid of all column monomial matrices. It is clear that as in the section on the Dowling lattice, by considering $G \wr S_{X}$ to be the monomial $|X| \times |X|$ matrices over $G$, this groups acts by conjugation on $B(G,X)$. It is clear that we can extend this action to the collection of all trivial subgroupoids of $B(G,X)$ which is an invariant set. Therefore, we have a natural action of $G \wr S_{X}$ on $\wh{R}_{X}(G)$ by the previous theorem.

\section{Review of the theory of BRSC}
\label{dar}

We will use the theory of Boolean representable simplicial complexes in order to study the matroids associated to both the Dowling and Rhodes lattices. We collect the basic facts in this section. See \cite{brsc} for more details.

Given a set $V$ and $n \geq 0$, we denote by $P_n(V)$
(respectively $P_{\leq n}(V)$) the set
of all subsets of $V$ with precisely (respectively at most) $n$
elements. To simplify notation, we shall often represent sets $\{ a_1, a_2,
\ldots, a_n\}$ in the form $a_1a_2\ldots a_n$.

A (finite) simplicial complex is a structure of the form $\H\; =
(V,H)$, where $V$ is a finite nonempty set and $H \subseteq 2^V$ is
nonempty and closed under taking subsets. Simplicial complexes, in
this abstract viewpoint, are also known as {\em hereditary collections}. The subsets in $H$ are said to be {\em independent}. A maximal independent set is called a {\em basis}. The maximum size of a basis is the {\em rank} of $\H$. A minimal dependent set id called a {\em circuit}.

Two simplicial complexes $(V,H)$ and $(V',H')$ are {\em isomorphic} if
there exists a bijection $\p:V \to V'$ such that
$$X \in H \mbox{ if and only if }X\p \in H'$$
holds for every $X \subseteq V$.

If $\H\; = (V,H)$ is a simplicial complex and $W \subseteq V$ is
nonempty, we call
$$\H|_W = (W,H \cap 2^W)$$
the {\em restriction} of $\H$ to $W$. It is obvious that $\H|_W$ is
still a simplicial complex.

We say that $X
\subseteq V$ is a {\em flat} of $\H$ if
$$\forall I \in H \cap 2^X \hspace{.2cm} \forall p \in V \setminus X
\hspace{.5cm} I \cup \{ p \} \in H.$$
The set of all flats of $\H$ is denoted by
$\flatx\H$.

Clearly, the intersection of any set of flats (including $V =
\cap\emptyset$) is still a flat. If we order $\flatx\H$ by inclusion,
it is then a $\wedge$-semilattice, and therefore a lattice with
$$(X \vee Y) = \cap\{ F \in \flatx\H \mid X \cup Y \subseteq F \}$$
for all $X,Y \in \flatx\H$.
We call $\flatx\H$ the {\em lattice of flats} of $\H$.
The lattice of flats induces a closure operator on $2^V$ defined by
$$\oo{X} = \cap\{ F \in \flatx \H \mid X \subseteq F \}$$
for every $X \subseteq V$.

We say that $X$ is a {\em transversal of the
successive differences} for a chain of subsets of $V$ of the form
$$A_0 \subset A_1 \subset \ldots \subset A_k$$
if $X$ admits an enumeration $x_1,\ldots , x_k$ such that $x_i \in A_i
\setminus A_{i-1}$ for $i = 1,\ldots,k$.

Let $\H\; = (V,H)$ be a simplicial complex. If $X \subseteq V$ is a
{\em transversal of the successive differences} for a chain
$$F_0 \subset F_1 \subset \ldots \subset F_k$$
in $\flatx \H$, it follows easily by induction that $\{ x_1, x_2,
\ldots, x_i \} \in H$ for $i = 0,\ldots,k$. In particular, $X \in H$.

We say that $\H$ is {\em
  Boolean representable} if every $X \in H$ is a transversal of the
successive differences for a chain in $\flatx \H$. We denote by $\BoR$
the class of all (finite) Boolean representable simplicial
complexes (BRSC). It is proved in \cite[Theorem 3.6.2]{brsc} that if $\H = (V,H)$
is a BRSC, then a subset
$X$ of $V$ is in $H$ if and only if $X$ admits some enumeration $X = x_{1},\ldots, x_{k}$
such that $B < \oo{x_{1}} < (\oo{x_{1}} \vee \oo{x_{2}}) < \ldots < (\oo{x_{1}} \vee \ldots \vee \oo{x_{k}})$ is a strict
chain in $\flatx \H$, where $B$ is the closure of the empty set.

A simplicial complex $\H\; = (V,H)$ is called a {\em matroid} if it
satisfies the {\em exchange property}:
\bi
\item[(EP)]
For all $I,J \in H$ with $|I| = |J|+1$, there exists some
  $i \in I\setminus J$ such that $J \cup \{ i \} \in H$.
\ei

A matroid is {\em simple} if every two element set is independent. Let $M$ be a matroid on a vertex set $V$. It is well known that by identifying two element sets that are not independent, we obtain a simple matroid $\widehat{M}$ whose lattice of flats $L(\widehat{M})$ is isomorphic to the lattice of flats of $M$ \cite{Oxley}. Furthermore, the lattice of flats of a matroid $M$ is a {\em geometric lattice}. Geometric lattices are lattices that are both atomistic and semimodular. We recall the definitions of these properties.

An atom $v$ of a lattice $L$ is a cover of the bottom element of $L$. Let $At(L)$ denote the set of atoms of a lattice $L$. A lattice is {\em atomistic} if every element $v \in L$ is a join of atoms, which are necessarily precisely the atoms below $v$. That is, for all $v \in L$, $v = \bigvee a$, where the join is over all atoms $a$ below $v$ in $L$. A lattice $L$ is {\em semimodular} if, whenever $v,w \in L$ and $v$ covers the meet $v \wedge w$, then the join $v \vee w$ covers $w$. A finite lattice $L$ is {\em graded}, if any two maximal chains between two given elements have the same length. If $L$ is a finite lattice we define the rank function $r:L \rightarrow \mathbb{N}$ by letting $r(v)$ be the length of the longest chain from the bottom of the lattice to $v$. It is known that a finite lattice is semimodular if and only if it is graded and the rank function satisfies
$r(v)+ r(w) \geq r(v \vee w) + r(v \wedge w)$ for all $v,w \in L$.

A classic theorem of Birkhoff \cite{Oxley} allows one to recover a simple matroid $M$ from its (geometric) lattice of flats $L(M)$. That is, every geometric lattice is the lattice of flats of a unique simple matroid. Indeed, let $L$ be a geometric lattice. Let $V=At(L)$ be the set of atoms of $L$. We define a simplicial complex $M(L)$ over $V$ by: $I \subseteq V$ is independent if and only if $r(\bigvee_{v \in I}v)=|I|$. Then Birkhoff's Theorem says that $M(L)$ is a simple matroid. Furthermore, the lattice of flats of $M(L)$ is isomorphic to $L$ and if $M$ is a simple matroid, and $L$ is the lattice of flats of $M$, then $M$ is isomorphic to $M(L)$.

Assuming that $L$ is geometric, let $I =\{e_{1}, \ldots, e_{k}\}$ be a subset of $V$. It is clear that $I$ is independent in $M(L)$ if and only if
$B < e_{1}< (e_{1} \vee e_{2}) < \ldots < (e_{1} \vee \ldots \vee e_{k})$ is a chain, where $B$ is the bottom element of $L$. In this case the ordering $(e_{1},e_{2} \ldots, e_{k})$ provides a transversal of the successive differences of this chain. Conversely, if $I =\{e_{1}, \ldots e_{k}\} \subseteq V$ is a transversal of the successive differences of any chain in $L$, then $r((e_{1} \vee \ldots \vee e_{k}))=k$ and $I$ is in $M(L)$. It follows that every simple matroid is Boolean representable. It can be shown that in fact every matroid is Boolean representable.

One of the important innovations of the theory of BRSC is that we can find lattices smaller than the lattice of flats that represent a simplicial complex. Let $L$ be an arbitrary finite lattice and let $V$ be a subset of $L$ that generates $L$ with respect to join. That is, every element of $L$ is a join of elements from $V$. We can define the analogue of a transversal of successive differences in $L$ with respect to $V$. This defines a simplicial complex $S(L,V)$, where the independent sets are precisely the subsets of $V$ that have some ordering which is a transversal of the successive differences. Equivalently \cite[Theorem 3.6.2]{brsc}, $S(L,V)$ consists of the set of subsets $X$ of $V$ that admit some enumeration that gives a strict chain in $L$. It is known that $S(L,V)$ is Boolean representable and thus the collection of all $S(L,V)$ coincides with the collection of BRSC \cite{brsc}. A lattice representation of a BRSC $S$ is a pair $(L,V)$ of a lattice $L$ with a join generating set $V$ such that $S=S(L,V)$. In a precise sense, the representation of a BRSC by its lattice of flats is its largest representation, containing all other representations \cite{brsc}. This allows for representations of simplicial complexes with lattices whose order is much smaller than that of the lattice of flats of the complex.

\be

Let $L_n$ be the chain $ 1 < 2 < \ldots < n$ and let $V_n = \{1, \ldots, n\}$. Then $S(L_n,V_n)$ = $U_{n,n}$, the uniform matroid where all subsets are independent.

\ee

It is clear that the lattice of flats of $U_{n,n}$ is the lattice of all subsets of an $n$ element set. We see from this example, that a matroid can be represented by a non-geometric lattice and that the size of a representing lattice can be exponentially smaller than the size of the lattice of flats.

\section{Frame matroids and lift matroids}
\label{frali}

For all concepts and results in this section, the reader is referred to \cite{Zaslavsky2}.
Throughout this section, all (undirected) graphs are allowed loops and multiple edges. A {\em theta} is the union of two cycles whose intersection is a nontrivial path:

\bigskip

$$
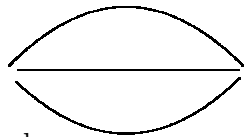
\xymatrix{
\ \ar@{-}@/^2.0pc/[rrr] \ar@{-}@/_2.0pc/[rrr] \ar@{-}[rrr] &&&
}$$ 

\bigskip

\noindent
Thus a theta possesses precisely three cycles. 

We say that a connected graph $\Gamma$ is {\em unicyclic} if the number $v$ of vertices equals the number $e$ of edges.
Since finite trees can be characterized as connected graphs satisfying $e = v-1$, 
it follows that unicyclic graphs are precisely those graphs which can be obtained by adding a new edge to a tree. 

A {\em biased graph} is a structure of the form $(\Gamma, \C)$, where $\Gamma$ is a (finite undirected) graph and $\C$ is a collection of cycles of $\Gamma$ satisfying the {\em theta property}: whenever two cycles in a theta belong to $\C$, so does the third.

A cycle is {\em balanced} if it belongs to $\C$. Otherwise, it is {\em unbalanced}. A subgraph is unbalanced if it contains an unbalanced cycle, and is fully unbalanced if all its cycles are unbalanced.

Given a biased graph $(\Gamma,\C)$, with edge set $E$, we define the {\em frame matroid} $F(\Gamma,\C) = (E,H)$ as follows: given $X \subseteq E$, we have $X \in H$ if and only if each connected component of the subgraph defined by $X$ is either a tree or an unbalanced unicyclic graph.

It is easy to check that $X \subseteq E$ is a circuit of $F(\Gamma,\C)$ if and only if $X$ is one of the following:
\bi
\item
a balanced cycle,
\item
the union of two unbalanced cycles sharing a vertex,
\item
the union of two vertex ­disjoint unbalanced cycles with a minimal path joining them,
\item
a fully unbalanced theta.
\ei
Graphs of the second and third types
$$
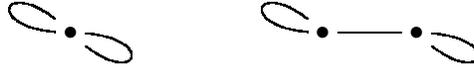
\xymatrix{
\bullet \ar@{-}@(ul,l) \ar@{-}@(dr,r) &&&
\bullet \ar@{-}@(ul,l) \ar@{-}[r] & \bullet \ar@{-}@(dr,r) 
}$$ 
are known as {\em tight} and {\em loose handcuffs}, respectively.

On the other hand, $X \subseteq E$ is a flat of $F(\Gamma,\C)$ if and only if $X$ satisfies the following conditions:
\bi
\item
if $C \in \C$, then $|C\setminus X| \neq 1$;
\item
if $C$ is an unbalanced cycle of $\Gamma$ and $Y$ is an unbalanced connected component of $X$, then $|C\setminus Y| \neq 1$;
\item
no two unbalanced connected components of $X$ are connected by an edge in $\Gamma$.
\ei

Next we define the {\em lift matroid} $L(\Gamma,\C) = (E,H')$: given $X \subseteq E$, we have $X \in H$ if and only if each connected component of the subgraph defined by $X$ is either a tree or an unbalanced unicyclic graph, and there is at most one component of the second type.

It is easy to check that $X \subseteq E$ is a circuit of $L(\Gamma,\C)$ if and only if $X$ is one of the following:
\bi
\item
a balanced cycle,
\item
the union of two unbalanced cycles sharing a vertex,
\item
the disjoint union of two unbalanced cycles;
\item
a fully unbalanced theta.
\ei
On the other hand, $X \subseteq E$ is a flat of $L(\Gamma,\C)$ if and only if $X$ satisfies the following conditions:
\bi
\item
if $C \in \C$, then $|C\setminus X| \neq 1$;
\item
if $C$ is an unbalanced cycle of $\Gamma$ and $X$ is unbalanced, then $|C\setminus X| \neq 1$.
\ei

An important example of biased graphs is given by gain graphs. Given a finite graph $\Gamma$ and a group $G$, we can construct a {\em gain graph} by associating elements of $G$ to the edges of $\Gamma$ with the help of an orientation. More precisely, given an edge $p \edge q$ in $\Gamma$, we associate a label $g \in G$ to the directed edge $p \mapright{} q$, and in this case we label the opposite edge $q \mapright{} p$ by $g\inv$. Now the label of a (directed) cycle 
$$p_1 \mapright{g_1} p_2 \mapright{g_2} \cdots \edge p_{m} \mapright{g_m} p_1$$
is $g_1\ldots g_m \in G$. The label of the cycle is well defined up to conjugacy and inversion. In particular, the label being 1 does not depend on neither of these factors. 

For such a gain graph, we define as balanced those cycles which have label 1.

\section{The relationship between the Rhodes and Dowling lattices for the trivial group and their matroids}

The purpose of this section is to clarify the relationship between the Dowling lattice $Q_{n}(1)$ and the Rhodes lattice $\wh{R}_{n}(1) = R_n(1)$. It is well known that the Dowling lattice $Q_{n}(1)$ is isomorphic to the lattice $\Pi_{n+1}$ of all (full) partitions of $\{0, \ldots n\}$. Let $\pi$ be a partial partition of $\{1, \ldots n\}$ such that the union of the classes is a subset $X$ of $\{1, \ldots n\}$. It is straightforward to see that the map that sends $\pi$ to the full partition  whose classes are those of $\pi$ together with $\{0,1,\ldots,n\} \setminus X$ is an isomorphism of lattices between $Q_{n}(1)$ and $\Pi_{n+1}$. In particular, $Q_{n}(1)$ is a geometric lattice and its corresponding (simple) matroid is the graphic matroid $\Gamma(K_{n+1})$ of the complete graph $K_{n+1}$ \cite{Oxley}.

On the other hand, the Rhodes lattice $R_{n}(1)$ is not even an atomistic lattice for $n \geq 2$. One needs only note that $(\{1,2\},\{\{1,2\}\})$ is a join irreducible element, as it covers uniquely the element $(\{1,2\},\{\{1\}, \{2\}\})$ which is the join of the two atoms $(\{1\},\{\{1\}\})$ and $(\{2\},\{\{2\}\})$. Despite this, $R_{n}(1)$ defines a matroid, which is the direct sum of the uniform matroid $U_{n,n}$ and the graphic matroid of the complete graph $K_n$. More precisely, $R_n({1})$ is a lattice representation of this matroid in the sense of the theory of BRSC \cite{brsc}.  

We will describe the BRSC represented by the Dowling lattices $Q_{n}(1)$ and the Rhodes lattices $R_{n}(1)$ with respect to their (unique) minimal sets of join generators, namely the join irreducible elements.

We have mentioned above that $Q_{n}(1)$ is isomorphic to the full partition lattice $P_{n+1}$. It is well known that $P_{n+1}$ is a geometric lattice and that its simple matroid is the graphic matroid of the complete graph $K_{n+1}$. This is the matroid whose vertices are the edges of $K_{n+1}$ and whose independent sets are the forests.

We mentioned above that $R_{n}(1)$ is not atomistic and thus not a geometric lattice if $n \geq 2$. Nonetheless we will see that the simplicial complex represented by $R_{n}(1)$ is a matroid. We first describe the join irreducible elements of $R_{n}(1)$.

\begin{Lemma}\label{JI}

Let $A_n =\{(\{i\},\{\{i\}\}) \mid i = 1, \ldots n\} \cup \{(\{i,j\}, \{\{i,j\}\}) \mid 1 \leq i < j \leq n\}$. Then $A_n$ is the set of join irreducible elements of $R_{n}(1)$.

\end{Lemma}

\proof It is clear that $\{(\{i\},\{\{i\}\})|i = 1, \ldots n\}$ is the set of atoms in $R_{n}(1)$ and that the unique element covered by $(\{i,j\},\{\{i,j\}\})$ is
$(\{i,j\},\{\{i\},\{j\}\})$. Therefore $A_n$ is contained in the set of join irreducible elements of $R_{n}(1)$.

Let $(X,\pi) \in R_{n}(1) \setminus A_n$. If $|X|$ = 2, then $(X,\pi) = (\{i,j\},\{\{i\},\{j\}\}) = ((\{i\},\{\{i\}\}) \vee (\{j\},\{\{j\}\}))$ and $(X, \pi)$ is not join irreducible. So we can assume that $|X| > 2$. If $\pi$ has only singleton classes, then $(X,\pi) = \bigvee_{i \in X}(\{i\},\{\{i\}\})$ and $(X,\pi)$ is not
join irreducible. If $\pi$ is an atom in the partition lattice of $X$, then there are distinct $i,j \in X$ such that
$\pi = \{\{i,j\},\{\{k\} \mid k \in X \setminus \{i,j\}\}\}$ and there is at least one singleton class, since $X$ has at least 3 elements. Therefore,
$(X,\pi) = (\{i,j\},\{\{i,j\}\}) \vee (\bigvee_{x \in X \setminus \{i,j\}} (\{x\},\{\{x\}\}))$. Finally, if $\pi$ is not an atom in the partition lattice of $X$, then $(X,\pi)=\bigvee (X,a)$, where the join is over the collection of atoms of the partition lattice of $X$ whose join is $\pi$. This completes the proof.
\qed

We now describe the cover relation in $R_{n}(1)$. Recall that a cover in a poset $(P,<)$ is a pair $x,y \in P$ such that $x>y$ and for all $z \in P$, if
$x \geq z \geq y$, then either $x=z$ or $y=z$. The next lemma follows easily from the definitions and the proof is left to the reader.

\begin{Lemma} \label{Cover}

Let $(X,\pi) \in R_{n}(1)$. Then $(Y,\tau)$ covers $(X,\pi)$ if and only if one of the following conditions holds:

\begin{itemize}
\item[(i)] $Y =X \cup \{i\}$ for some $i \notin X$ and $\tau = \pi \cup\{\{i\}\}$ In this case, $(\{i\},\{\{i\}\})$ is the unique join irreducible element $a$ of $R_{n}(1)$ such that
$(Y,\tau) = ((X,\pi) \vee a)$.
\item[(ii)] $Y=X$ and $\tau$ covers $\pi$ in the partition lattice of $X$, that is two classes of $\pi$ are merged into a single class of $\tau$. In this case the set of join irreducible elements $a$ of $R_{n}(1)$ such that $(Y,\tau) = ((X,\pi) \vee a)$ is equal to the set of all $(\{i,j\},\{\{i,j\}\})$ such that $i$ and $j$ are in the distinct classes of $\pi$ that produce $\tau$.
\end{itemize}

\end{Lemma}

Thus covers either are a cover in the subset lattice $(2^n,\cup)$ or in the partition lattice $(\Pi_{n}, \leq)$. We can make this more precise by computing the BRSC defined by $R_{n}(1)$. Recall that if $S = (V,H)$ and $S' = (V',H')$ are simplicial complexes such that $V$ and $V'$ are disjoint then their direct sum is the simplicial complex $S \bigoplus S'$ with vertices $V \cup V'$ and independent subsets of the form $\{X \cup Y\mid X \in H, Y \in H'\}$. Both $(2^n,\cup)$ and $(\Pi_{n}, \leq)$ are geometric lattices. The BRSC of  $(2^n,\cup)$ is the matroid $U_{n,n}$, the uniform matroid of all subsets of $\underline{n}$. The BRSC of $(\Pi_{n}, \leq)$ is the graphic matroid $\Gamma(K_{n})$ of the complete graph $K_n$. This has the edges of $K_n$ as vertices and the forests of $K_n$ as independent sets \cite{Oxley}.

\bt

The BRSC of $R_{n}(1)$ with respect to its set $A_n$ of join irreducible elements is isomorphic to the direct sum $U_{n,n}\bigoplus \Gamma(K_{n})$.

\et

\proof Let $A_n$, the set of join irreducible elements of $R_{n}(1)$, be as in Lemma \ref{JI}. We first note that under the function that sends $(\{i\},\{\{i\}\})$ to
$\{i\}$ ($i=1 \ldots n$) and $(\{i,j\},\{\{i,j\}\})$ to the partition $\{\{i,j\},\{k\} \mid k \neq i,j\}$ is a bijection between $A_n$ and the disjoint union of
the sets of atoms of $(2^n,\cup)$ and $(\Pi_{n}, \leq)$. We will use this identification in the rest of the proof.

We next note that if $X = (\underline{n}, id)$ where $id$ is the partition consisting of all the singletons of $\underline{n}$, then the down ideal $X^\downarrow$ of all elements of $R_{n}(1)$ less than or equal to $X$ is a lattice isomorphic to $(2^n,\cup)$ and that the up ideal $X^\uparrow$ is a lattice isomorphic to $(\Pi_{n}, \leq)$.

Let $B$ be a basis element of $\Gamma(K_{n})$. Then $B$ corresponds to a spanning tree of $K_n$ and when translating to the language of $c$-independent sets in
$(\Pi_{n}, \leq)$ to a set of $n-1$ atoms of $(\Pi_{n}, \leq)$. The atom $\{\{i,j\},\{k\} \mid k \neq i,j\}$ corresponds to the edge $\{i,j\}$. We claim that the set (under the identification in the last paragraph), $\{\{i\}\mid i=1, \ldots n\} \cup B$ is $c$-independent in $R_{n}(1)$. Indeed, $\{\{i\}\mid i=1, \ldots n\}$ labels a chain from the bottom of $R_{n}(1)$,  to $X = (\underline{n}, id)$ and $B$ labels a chain from $X$ to the top of $R_{n}(1)$, again for any enumeration. This is clearly a basis element of the BRSC of $R_{n}(1)$ and thus, every basis element of $U_{n,n}\bigoplus \Gamma(K_{n})$ is also a basis element of $R_{n}(1)$, under the identification of the previous paragraph.

Conversely, let $Y=\{x_{i}\mid i=1 \ldots k\}$ be a set of join irreducible elements that is a basis for the BRSC of $R_{n}(1)$. We can assume that the enumeration $x_{1}, \ldots, x_{k}$ defines a chain, $x_{1} < \ldots < (x_{1} \vee \ldots \vee x_{k})$ in $R_{n}(1)$. Since $Y$ is a basis, each inequality in this chain must be a cover and
$(x_{1} \vee \ldots \vee x_{k})$ is the top element of $R_{n}(1)$. Since the only join irreducible elements that enlarge the set component in $R_{n}(1)$ are atoms of the form
($\{i\},\{\{i\}\})$, it follows from Lemma \ref{Cover} that $Y$ must contain all such atoms, $i=1, \ldots n$. Similarly, since the only join irreducible elements that enlarge the partition
component without enlarging the set component are of the form, $(\{i,j\},\{\{i,j\}\})$, it follows that $Y$ must contain a basis of the matroid of  $\Pi_{n}$ under the identification above. It
follows that under this identification, $Y$ corresponds to a basis of $U_{n,n}\bigoplus \Gamma(K_{n})$. Therefore,  this identification induces a bijection between the
bases of the BRSC of $R_{n}(1)$ and $U_{n,n}\bigoplus \Gamma(K_{n})$. Therefore, these two BRSC are isomorphic.
\qed

\section{The matroid defined by $R_n(G)$}
\label{vlk}

%

We consider now $R_n(G) = R_X(G)$ for $X = \un{n}$. Let $\omega_I$ (respectively $\iota_I$) denote the partition of a set $I$ with one single block(respectively, only singleton blocks). We omit the subscript when the set is implicit.
Write 
\bi
\item
$B_n = \{ (\{ i \}, \omega,[f]_{\omega}) \mid i = 1,\ldots,n;\, f \in F(\{ i\},G) \}$, 
\item
$C_n = \{ (\{ i,j\}, \omega,[f]_{\omega}) \mid 1 \leq i < j \leq n;\,  f \in F(\{ i,j\},G) \}$,
\item
$A_n = B_n \cup C_n$. 
\ei

A straightforward adaptation of Lemma \ref{JI} yields:

\bl
\label{sjig}
$A_n$ is the set of join irreducible elements of $R_n(G)$.
\el

Since the join irreducible elements of $\wh{R}_n(G)$ are precisely the join irreducible elements of $R_n(G)$, we also get:

\bc
\label{sjigl}
Let $n,|G| > 1$. Then $A_n$ is the set of
 join irreducible elements of $\wh{R}_n(G)$.
\ec

Let $\H_n(G) = (A_n,H_n(G))$ denote the simplicial complex $\H_n(G) = (A_n,H_n(G))$ as follows: $H_n(G)$ is the set of all subsets of $A_n$ admitting an enumeration $\alpha_1,\ldots,\alpha_m$ such that
$$\alpha_1 <_R (\alpha_1 \vee \alpha_2) <_R \ldots <_R (\alpha_1 \vee \ldots \vee \alpha_m).$$
Note that this implies that $\alpha_1,\ldots, \alpha_m$ have a common upper bound, which implies that joins are well defined throughout the chain. This is a well-defined simplicial complex since $\emptyset \in H_n(G)$ and $H_n(G)$ is closed under taking subsets. This is of course reminiscent of the BRSC defined in Section \ref{dar} from an arbitrary lattice. Although $R_n(G)$ is in general just a semilattice, $\H_n(G)$ turns out to be Boolean representable as well. Indeed, a matroid.

We aim at describing the bases of $\H_n(G)$. It is useful to recall the proof of Theorem \ref{tsg}, which allows us to view the elements of $R_n(G)$ as trivial subgroupoids of $B(G,\un{n})$ when convenient. Trivial subgroupoids may be pictured as labeled directed graphs, where the (single) morphism between objects $i,j \in \un{n}$ (if it exists), is pictured by $i \mapright{g} j$ for some $g \in G$. In view of the rules of composition of morphisms in groupoids, we must have a morphism $i \mapright{gg'} k$ whenever $i \mapright{g} j$ and $j \mapright{g'} k$ are morphisms. It follows from uniqueness that $g = 1$ necessarily if $i = j$.
Moreover, since a morphism $i \mapright{g} j$ is supposed to be an isomorphism, and $gg' = 1$ whenever $i \mapright{gg'} i$ is a morphism, we have also a morphism $j \mapright{g\inv} i$. This perspective will allow us to use graph-theoretical tools to describe independent sets and bases.

In the subgroupoid perspective, we need to identify the elements of $B_n$ and $C_n$, which are described respectively by
$$\xymatrix{
i \ar@(ul,l)_1 &
(i \in \un{n})&\mbox{and}&i \ar@/^/[rr]^g&& j \ar@/^/[ll]^{g\inv}&(1 \leq i < j \leq n,\, g \in G)
}$$
Assuming the necessity of the inverse edges, it is of course enough to represent the above element of $C_n$ by the one-edge graph $i \mapright{g} j$. Now each $Z \subseteq C_n$ can be represented by a directed labeled graph $\Gamma(Z)$. We denote by $\Gamma_0(Z)$ the undirected (multi)graph with the same vertex set as $\Gamma(Z)$ and an edge $i \edge j$ for each pair of inverse edges
$$\xymatrix{
i \ar@/^/[rr]^{g}&& j \ar@/^/[ll]^{g\inv}
}$$
in $\Gamma(Z)$.

\be
Let $Z = \{ (\{ 1,2\}, \omega,[f]_{\omega}), (\{ 2,3\}, \omega,[h]_{\omega}), (\{ 2,3\}, \omega,[h']_{\omega}) \}$, where $f(1) = f(2) = h(2) = h(3) = h'(3) = 1$ and $h'(2) = g$. Then $\Gamma(Z)$ is the graph
$$\xymatrix{
1 \ar[rr]^1 && 2 \ar@/^/[rr]^{(1)} \ar@/_/[rr]_{g\inv} && 3 
}$$
and  $\Gamma_0(Z)$ is the graph
$$\xymatrix{
1 \ar@{-}[rr] && 2 \ar@{=}[rr] && 3 
}$$
\ee

Let 
$$\H^{(1)}_n(G) = (B_n,H^{(1)}n(G)) = \H_n(G)|_{B_n} \quad \mbox{ and } \H^{(2)}_n(G) = (C_n,H^{(2)}_n(G)) = \H_n(G)|_{C_n}.$$

\bl
\label{dsum}
\bi
\item[(i)] $\H_n(G) \cong \H^{(1)}_n(G) \bigoplus \H^{(2)}_n(G)$.
\item[(ii)] $\H^{(1)}_n(G) \cong U_{n,n}$.
\ei
\el

\proof
(ii) Since $|B_n| = n$, it suffices to show that $B_n \in H^{(1)}n(G)$. For each $i \in \un{n}$, write $\alpha_i = (\{ i \}, \omega,[f_i]_{\omega})$ with $f_i(i) = 1$. It is immediate that
\beq
\label{dsum1}
(\alpha_1 \vee \ldots \vee \alpha_i) = (\{ 1,\ldots, i \}, \iota, [f_1 \cup \ldots \cup f_i]_{\iota}),
\eeq
whence 
$$\alpha_1 <_R (\alpha_1 \vee \alpha_2) <_R \ldots <_R (\alpha_1 \vee \ldots \vee \alpha_n).$$
Thus $B_n \in H^{(1)}_n(G)$ and so $\H^{(1)}n(G) \cong U_{n,n}$.

(i) Write
$$H' = \{ Y \cup Z \mid Y \in H^{(1)}n(G),\, Z \in H^{(2)}n(G) \}.$$
We need to show that $H_n(G) = H'$.
Let $X \in H_n(G)$. Then
$X = (X \cap B_n) \cup (X \cap C_n)$. Since $X\cap B_n \in H^{(1)}n(G)$ and $X\cap C_n \in H^{(2)}_n(G)$, it follows that $X \in H'$. Thus $H_n(G) \subseteq H'$.

Conversely, let $X \in H'$. Then $X = Y \cup Z$ for some $Y \in H^{(1)}_n(G)$ and $Z \in H^{(2)}_n(G)$. We may assume by (ii) that $Y = B_n = \{ \alpha_1,\ldots,\alpha_n\}$.
Assume that $Z = \{ \beta_1,\ldots,\beta_m\}$ and 
$$\beta_1 <_R (\beta_1 \vee \beta_2) <_R \ldots <_R (\beta_1 \vee \ldots \vee \beta_m).$$
Let $\alpha = (\un{n}, \iota,[f_1 \cup \ldots \cup f_n]_{\iota})$. By (\ref{dsum1}), we have
$\alpha = (\alpha_1 \vee \ldots \vee \alpha_n)$. We claim that
\beq
\label{dsum2}
(\alpha \vee \beta_1) <_R (\alpha \vee \beta_1 \vee \beta_2) <_R \ldots <_R (\alpha \vee \beta_1 \vee \ldots \vee \beta_m).
\eeq
Together with (\ref{dsum1}), this implies that $X = Y \cup Z \in H_n(G)$, yielding $H' = H_n(G)$.

Let $i \in \{2,\ldots,m\}$. Assume that $(\beta_1 \vee \ldots \vee \beta_{i-1}) = (I,\pi,[f]_{\pi})$ and $(\beta_1 \vee \ldots \vee \beta_{i}) = (J,\tau,[h]_{\tau})$. Straightforward computation shows that 
$$((I,\pi,[f]_{\pi}) \vee \alpha) = (\un{n},\pi',[f']_{\pi'}),$$
where $\pi'$ is the partition of $\un{n}$ obtained by adding singleton classes to $\pi$, and $f'|_I = f$. Similarly, we may write
$$((J,\tau,[h]_{\tau}) \vee \alpha) = (\un{n},\tau',[h']_{\tau'}).$$
Clearly, $(I,\pi,[f]_{\pi}) <_R (J,\tau,[h]_{\tau})$ yields 
$$(\un{n},\pi',[f']_{\pi'}) = ((I,\pi,[f]_{\pi}) \vee \alpha) \leq_R ((J,\tau,[h]_{\tau}) \vee \alpha) = (\un{n},\tau',[h']_{\tau'}).$$
Suppose that 
$(\un{n},\pi',[f']_{\pi'}) = (\un{n},\tau',[h']_{\tau'})$. Then $\pi' = \tau'$. Note that, in a join of elements from $C_n$, there are no singleton classes. Therefore, by removing the singleton classes from $\pi'$ (respectively $\tau'$), we obtain precisely $\pi$ (respectively $\tau$). Thus $\pi' = \tau'$ yields $\pi = \tau$. It follows that also $I = J$ and $[f]_{\pi} = [h]_{\tau}$, contradicting $(I,\pi,[f]_{\pi}) <_R (J,\tau,[h]_{\tau})$.

Therefore $(\un{n},\pi',[f']_{\pi'}) <_R (\un{n},\tau',[h']_{\tau'})$ and so (\ref{dsum2}) holds as required.
\qed

Before discussing the independent sets of $\H^{(2)}_n(G)$, we prove the following lemma.

\bl
\label{thefaces1}
Let $Z \subseteq C_n$. If $\Gamma_0(Z)$ is a forest with vertex set $I$, then $\vee Z = (I,\pi,[f]_{\pi})$ for some partition $\pi$ of $I$ and some $f \in F(I,G)$.
\el

\proof
We use induction on the number of edges of $\Gamma_0(Z)$. 
If $\Gamma_0(Z)$ has 0 edges, then $Z = \emptyset$ and so $\vee Z = (\emptyset,\emptyset,[\emptyset]_{\emptyset})$. Hence the lemma holds for 0 edges.

Assume now that $\Gamma_0(Z)$ has $m > 0$ edges and the lemma holds for forests with less edges. Let $j \in I$ have degree 1 in $\Gamma_0(Z)$. Then there exists a unique $i \in I$ and a unique $g \in G$ such that $i \mapright{g} j$ represents an element of $Z$, say $(\{i,j\}, \omega,[h]_{\omega})$ (so that $g = (h(i))\inv h(j)$). Let $Z'$ be obtained from $Z$ by removing this element. We can have one of two cases:

\medskip
\noindent
{\bf Case 1:} $i$ has degree 1 in $\Gamma_0(Z)$.

\smallskip
\noindent
Then $\Gamma_0(Z')$ is a forest with vertex set $I\setminus \{i,j\}$ and $m-1$ edges. By the induction hypothesis, $\vee Z' = (I\setminus \{i, j\},\pi',[f']_{\pi'})$ for some partition $\pi'$ of $I\setminus \{ i,j\}$ and some $f' \in F(I\setminus \{ i,j\},G)$. Writing $\tau = \pi' \cup \{ \{ i,j \}\}$, it is now straightforward to check that
$$\vee Z = ((I\setminus \{i, j\},\pi',[f']_{\pi'}) \vee (\{i,j\}, \omega,[h]_{\omega})) = (I,\tau, [f'\cup h]_{\tau})$$
and so the lemma holds in this case.

\medskip
\noindent
{\bf Case 2:} $i$ has degree $> 1$ in $\Gamma_0(Z)$.

\smallskip
\noindent
Then $\Gamma_0(Z')$ is a forest with vertex set $I\setminus \{j\}$ and $m-1$ edges. By the induction hypothesis, $\vee Z' = (I\setminus \{j\},\pi',[f']_{\pi'})$ for some partition $\pi'$ of $I\setminus \{ j\}$ and some $f' \in F(I\setminus \{j\},G)$. It is now straightforward to check that
$$\vee Z = ((I\setminus \{j\},\pi',[f']_{\pi'}) \vee (\{i,j\}, \omega,[h]_{\omega})) = (I,\pi, [f]_{\pi}),$$
where $\pi$ is the partition of $I$ obtained by adjoining $j$ to the block of $\pi'$ containing $i$, and $f$ extends $f'$ with $f(j) = (f'(i))g$. 
\qed

\bt
\label{thefaces}
$H^{(2)}_n(G) = \{ Z \subseteq C_n \mid \Gamma_0(Z) \mbox{ is a forest}\}.$
\et

\proof
We prove the opposite inclusion of the theorem by induction on the number of edges of $\Gamma_0(Z)$. 
If $\Gamma_0(Z)$ has 0 edges, then $Z = \emptyset$ and so $Z \in H^{(2)}_n(G)$ trivially. Assume now that $\Gamma_0(Z)$ has $m > 0$ edges and the claim holds for forests with $m-1$ edges. Let $j \in I$ have degree 1 in $\Gamma_0(Z)$. 
Then there exists a unique $i \in I$ and a unique $g \in G$ such that $i \mapright{g} j$ represents an element of $Z$, say $(\{i,j\}, \omega,[h]_{\omega})$ (so that $g = (h(i))\inv h(j)$). Let $Z'$ be obtained from $Z$ by removing this element. Since $\Gamma_0(Z)$ is a forest with vertex set contained in $I\setminus \{ j\}$, we have $Z' \in H^{(2)}_n(G)$ by the induction hypothesis. On the other hand,
it follows from Lemma \ref{thefaces1} that $\vee Z$ is of the form $(I,\ldots,\ldots)$, and $\vee Z'$ is of the form $(I',\ldots,\ldots)$ with $I' \subseteq I\setminus \{ j\}$. Together with $\vee Z' \leq_R \vee Z$, this implies $\vee Z' <_R \vee Z$. Since $Z' \in H^{(2)}_n(G)$, we get $Z \in H^{(2)}_n(G)$ and this completes the proof of the opposite inclusion.

Now let $Z \subseteq C_n$ and suppose that $\Gamma_0(Z)$ contains a cycle. Then $\Gamma(Z)$ contains a subgraph of the form
$$i_0 \mapright{g_1} i_1 \mapright{g_2} \ldots  \mapright{g_m} i_m = i_0$$
for some $m \geq 2$, where $(g_{j+1},i_{j+1}) \neq (g_j\inv,i_{j-1})$ for every $0 < j < m-1$. Write $I = \{ i_0,\ldots,i_{m-1} \}$ and let $Z' \subseteq Z$ correspond to the  edges of this cycle. 

Suppose that $Z' \in H^{(2)}_n(G)$. Then there exists an enumeration $\beta_1,\ldots,\beta_m$ of the elements of $Z'$ such that 
\beq
\label{thefaces2}
\beta_1 <_R (\beta_1 \vee \beta_2) <_R \ldots <_R (\beta_1 \vee \ldots \vee \beta_m).
\eeq
Without loss of generality, we may assume that $\beta_m$ corresponds to the edge $i_{m-1} \mapright{g_m} i_0$. It follows that $\beta_m = (\{ i_{m-1},i_0\}, \omega, [h]_{\omega})$ for some $h \in F(\{ i_{m-1},i_0\},G)$.

Let $Z''$ be obtained by removing from $Z'$ the element corresponding to this edge. Since $\Gamma_0(Z'')$ is a tree with vertex set $I$, it follows from Lemma \ref{thefaces1} that $\vee Z'' = (I,\pi,[f]_{\pi})$ for some partition $\pi$ of $I$ and some $f \in F(I,G$.  In fact, since each edge $i_{j-1} \edge i_j$ in $\Gamma_0(Z'')$ forces $i_{j-1}$ and $i_j$ to end up in the same class of $\pi$, it follows that $\pi = \{ I \}$. On the other hand, $(f(i_{j-1}))\inv f(i_j) = g_j$ for $j = 1,\ldots,m-1$. Hence $(f(i_0))\inv f(i_{m-1}) = g_1g_2\ldots g_{m-1}$.

We can have one of two cases:

\medskip
\noindent
{\bf Case 1:} $g_1g_2\ldots g_m = 1$.

\smallskip
\noindent
We show that 
$$\beta_m = (\{ i_{m-1},i_0\}, \omega,[h]_{\omega}) \leq_R (I,\omega,[f]_{\omega}) = \vee Z''.$$
Indeed, $\{ i_{m-1},i_0\} \subseteq I$ and 
$$(f(i_{m-1}))\inv f(i_0) = g_{m-1}\inv\ldots g_1\inv = g_m = (h(i_{m-1}))\inv h(i_0),$$
hence $[f|_{\{ i_{m-1},i_0\} }]_{\omega} = [h]_{\omega}$. This contradicts (\ref{thefaces2}), thus $Z' \notin H^{(2)}_n(G)$ in this case.

\medskip
\noindent
{\bf Case 2:} $g_1g_2\ldots g_m \neq 1$.

\smallskip
\noindent
Let $\beta'_m = (\{ i_{m-1},i_0\}, \omega,[h']_{\omega})$, where $(h'(i_{m-1}))\inv h'(i_0) = g_{m-1}\inv\ldots g_1\inv$. It follows from the discussion of Case 1 that $\beta'_m \leq_R \vee Z''$. By (\ref{thefaces2}), we get $\beta'_m \leq_R (\beta_1 \vee \ldots \vee \beta_m)$. But by Proposition \ref{existjoin} there is no common upper bound for $\beta_m$ and $\beta'_m$, thus we get a contradiction and so $Z' \notin H^{(2)}_n(G)$ in all cases. Since $Z' \subseteq Z$, we get $Z \notin H^{(2)}_n(G)$ as required.
\qed

We can now show that $\H_n(G)$ is a matroid and identify it by means of classical constructions, in fact it turns out to be the direct sum of a uniform matroid with the graphical matroid of a complete multigraph. Given $m,n \geq 1$, we denote by $mK_n$ the {\em complete multigraph} having $n$ vertices and $m$ (undirected) edges connecting any two distinct vertices.

\bc
\label{grama}
$\H_n(G) \cong U_{n,n} \bigoplus \Gamma(|G|K_{n})$.
\ec

\proof
In view of Lemma \ref{dsum}, it suffices to show that $\H^{(2)}_n(G) \cong \Gamma(|G|K_{n})$. Now $C_n$ can be seen as the set of edges of the complete multigraph $|G|K_n$, and the independent sets of the graphical matroid $\Gamma(|G|K_{n})$ are precisely the forests of $|G|K_n$. By Theorem \ref{thefaces}, these are precisely the independent sets of $\H^{(2)}_n(G)$.
\qed

Since matroids are obviously closed under direct sum, it follows that $\H_n(G)$ is a matroid. However, $\H_n(G)$ does not characterize $G$ since $\H_n(G) \cong \H_n(G')$ whenever $|G| = |G'|$.

\section{The matroid defined by $\wh{R}_n(G)$}
\label{alift}

In this section, we assume that $n,|G| > 1$. Whenever possible, we adapt results from Section \ref{vlk} with more concise explanations.

We consider now $\wh{R}_n(G) = \wh{R}_X(G)$ for $X = \un{n}$. By Corollary \ref{sjigl}, $A_n$ is the set of join irreducible elements of $\wh{R}_n(G)$. 
Let $\wh{\H}_n(G) = (A_n,\wh{H}_n(G))$ be the BRSC defined by $\wh{R}_n(G)$.

Let 
$$\wh{\H}^{(1)}_n(G) = (B_n,\wh{H}^{(1)}_n(G)) = \wh{\H}_n(G)|_{B_n} \quad \mbox{ and } \quad \wh{\H}^{(2)}_n(G) = (C_n,\wh{H}^{(2)}_n(G)) = \wh{\H}_n(G)|_{C_n}.$$

A straightforward adaptation of the proof of Lemma \ref{dsum} yields:
 
\bl
\label{dsumT}
\bi
\item[(i)] $\wh{\H}_n(G) \cong \wh{\H}^{(1)}_n(G) \bigoplus \wh{\H}^{(2)}_n(G)$.
\item[(ii)] $\wh{\H}^{(1)}_n(G) \cong U_{n,n}$.
\ei
\el


We denote by $\Delta_n(G)$ the gain graph obtained from the complete multigraph $|G|K_n$ by attributing all possible labels $g \in G$ to the $|G|$ distinct edges connecting each pair of distinct vertices.

Next we prove one of the main results of this section:

\bt
\label{lift}
$\wh{\H}^{(2)}_n(G)$ is the lift matroid of $\Delta_n(G)$.
\et

\proof
The mapping $\Gamma$ defined in Section \ref{vlk} establishes a bijection between $C_n$ and the edges of $|G|K_n$. The independent sets of a lift matroid were described in Section \ref{frali}, hence it suffices to show, for every $Z \subseteq C_n$, that 
$Z \in \wh{H}^{(2)}_n(G)$ if and only if the following conditions hold:
\bi
\item[(L1)]
each connected component of $\Gamma_0(Z)$ is either a tree or a unicyclic graph:
\item[(L2)]
at most one connected component of $\Gamma_0(Z)$ is unicyclic;
\item[(L3)]
if $\Gamma_0(Z)$ has a unicyclic connected component, then $\Gamma(Z)$ is unbalanced.
\ei

A chain of the form 
$$\alpha_1 <_R (\alpha_1 \vee \alpha_2) <_R \ldots <_R (\alpha_1 \vee \ldots \vee \alpha_m)$$
in $\wh{R}_n(G)$ is a chain in $R_n(G)$ unless $(\alpha_1 \vee \ldots \vee \alpha_m) = T$. Splitting 
$\wh{H}^{(2)}_n(G) = H^{(2)}_n(G) \cup (\wh{H}^{(2)}_n(G) \setminus H^{(2)}_n(G))$, we get two types for $Z \in \wh{H}^{(2)}_n(G)$:
\bi
\item
in view of Theorem \ref{thefaces}, each connected component of $\Gamma_0(Z)$ is a tree, in which case conditions (L1)--(L3) are verified;
\item
$Z$ is obtained by adding some $\alpha_{m+1} \in C_n$ to some $Z_0 = \{ \alpha_1, \ldots, \alpha_m \} \in H^{(2)}_n(G)$ so that $(\alpha_1 \vee \ldots \vee \alpha_m \vee \alpha_{m+1}) = T$. This means that $\Gamma_0(Z_0)$ is a forest but $\Gamma_0(Z)$ is not, and conditions (L1) and (L2) follow. By the arguments in Case 1 of the proof of Theorem \ref{thefaces}, the unique cycle in $\Gamma(Z)$ must be indeed unbalanced, otherwise the element corresponding to the last edge of the cycle lies below the join of the elements corresponding to the other edge. Thus (L3) must hold.
\ei
The converse implication is analogous.
\qed

\bc
\label{suli}
$\wh{\H}_n(G)$ is the direct sum of the complete matroid $U_{n,n}$ with the lift matroid of $\Delta_n(G)$.
\ec

\bc
\label{thefacetsT}
\bi
\item[(i)] $Z \subseteq C_n$ is a basis of $\wh{\H}^{(2)}_n(G)$ if and only if $\Gamma_0(Z)$ is a unicyclic graph with $n$ vertices and $\Gamma(Z)$ is unbalanced.
\item[(ii)] $X \subseteq A_n$ is a basis of $\wh{\H}_n(G)$ if and only if $B_n \subseteq X$ and $\Gamma_0(X\cap C_n)$ is a unicyclic graph with $n$ vertices and $\Gamma(Z)$ is unbalanced.
\item[(iii)] $\wh{\H}_n(G)$ has rank $2n$.
\ei
\ec

\proof
(i) By Theorem \ref{lift}, since any graph satisfying conditions (i) and (ii) can be extended to a unicyclic graph with vertex set $\oo{n}$.

(ii) By Lemma \ref{dsumT}(i), the bases of $\wh{\H}_n(G)$ are unions of a basis of $\wh{\H}^{(1)}_n(G)$ with a basis of $\wh{\H}^{(2)}_n(G)$. By Lemma \ref{dsumT}(ii), $A_n$ is the unique basis of $\wh{\H}^{(1)}_n(G)$.

(iii) A unicyclic graph with $n$ vertices has precisely $n$ edges. Thus each basis of $\wh{\H}_n(G)$ has $2n$ elements and therefore $\wh{\H}_n(G)$ has rank $2n$.
\qed

Given the extensive literature on lift matroids, we assume that the following results are known, but we include a proof for completeness

\bt
\label{countfacetsT}
The number of facets of $\wh{\H}_n(G)$  is precisely 
\beq
\label{countfacetsT2}
\left(
\frac{n!}{(n-2)!} n^{n-3} +\sum_{k=3}^n \frac{n!}{2(n-k)!} n^{n-k-1}\right)
|G|^{n-1}(|G|-1).
\eeq
\et

\proof
In view of Corollary \ref{thefacetsT}(ii), we have to count unicyclic graphs. In \cite[Section 6.8]{SF}, we can find a formula for the number of unicyclic graphs on $n$ vertices, but the authors exclude cycles of length 2:
\beq
\label{countfacetsT1}
\sum_{k=3}^n \frac{n!}{2(n-k)!} n^{n-k-1}.
\eeq
Adapting the arguments leading to this formula to the 2-cycle case, and using Cayley's formula for rooted forests \cite[Theorem 6.5]{SF}, we get that the number of unicyclic graphs on $n$ vertices is
$$\frac{n!}{(n-2)!} n^{n-3} +\sum_{k=3}^n \frac{n!}{2(n-k)!} n^{n-k-1}.$$
The absence of the factor 2 with respect to the summands of formula is easy to explain. The automorphism group of a $k$-cycle has order $2k$ if $k \geq 3$ and order $k$ if $k = 2$.

Now, if we consider a spanning tree for each one of these graphs (which has $n-1$ edges), we can choose freely the label in $G$, and for the remaining edge we have precisely $|G| -1$ choices (since we must ensure that the cycle becomes unbalanced). Therefore the number of bases of $\wh{\H}_n(G)$  is (\ref{countfacetsT2}).
\qed

We recall now a standard notation in asymptotics. Given mappings $f,g: \mathbb{N} \to \mathbb{R}$, we write $f(n) \sim g(n)$ if $\lim_{n\to +\infty} \frac{f(n)}{g(n)} = 1$.

\bc
\label{asyfT}
Let $\p(n)$ denote the number of bases of  $\wh{\H}_n(G)$. Then 
$$\p(n) \sim \sqrt{\frac{\pi}{8}} \, n^{n-\frac{1}{2}}|G|^{n-1}(|G|-1).$$
\ec

\proof
By \cite[Theorem 6.6]{SF}, we have $\sum_{k=3}^n \frac{n!}{2(n-k)!} n^{n-k-1} \sim \sqrt{\frac{\pi}{8}} \, n^{n-\frac{1}{2}}$. In view of Theorem \ref{countfacetsT}, it suffices to check that 
$$\lim_{n\to +\infty} \frac{\frac{n!}{(n-2)!} n^{n-3}}{\sqrt{\frac{\pi}{8}} \, n^{n-\frac{1}{2}}} = 0,$$
which is straightforward.
\qed

We can prove that, unlike $\H_n(G)$, the matroid $\wh{\H}_n(G)$ determines the group $G$ up to isomorphism.

\bt
\label{mdg}
Let $n,n' > 1$ and let $G,G'$ be finite nontrivial groups. Then the following conditions are equivalent:
\bi
\item[(i)] $\wh{R}_n(G) \cong \wh{R}_{n'}(G')$;
\item[(ii)] $\wh{\H}_n(G) \cong \wh{\H}_{n'}(G')$;
\item[(iii)] $n = n'$ and $G \cong G'$.
\ei
\et

\proof
(i) $\Rw$ (ii). Since $\wh{\H}_n(G)$ is fully determined by $\wh{R}_n(G)$ via its join irreducible elements.

(ii) $\Rw$ (iii). The matroid $\wh{\H}_n(G)$ determines $n$ by Corollary \ref{thefacetsT}(iii). It suffices to show that $\wh{\H}_n(G)$ determines $G$ up to isomorphism. Thus we assume that $\H = (V,H)$ is a matroid isomorphic to some $\wh{\H}_n(G)$, and try to reconstruct $G$ from $\H$. Note that $\wh{\H}_n(G)$ determines $|G|$ by Theorem \ref{countfacetsT}. Since nonisomorphic groups of the same order appear only at order 4, we may assume that $|G| \geq 4$.

We know that $V$ corresponds to the set $A_n$ of join irreducible elements of $\wh{\H}_n(G)$. It follows from Theorem \ref{lift} that any 3-subset of $A_n$ containing some element of $B_n$ belongs necessarily to $\wh{H}_n(G)$. Moreover, any element of $C_n$ is part of a balanced triangle which is not independent. Thus the subset of $V$ corresponding to $C_n$ may be identified as 
$$W = \{ w \in V \mid wvy \notin H \mbox{ for some distinct }v,y \in V \setminus \{ w\} \}.$$
Next we fix some element $a \in W$, which must correspond to some $i \mapright{g} j$ in $C_n$. 
The elements of $W$ corresponding to elements of the form $i \mapright{g'} j$ in $C_n$ are collected in
$$W_{ij} = \{ a \} \cup \{ w \in W \setminus \{ a \} \mid wav \notin H\mbox{ for at least two distinct }v \in W \setminus \{ a,w \} \}.$$
This follows from the fact that the unique elements of $P_3(C_n) \setminus \wh{H}_n(G)$ are of the form 
\beq
\label{vac}
i \mapright{g_1} j, \quad i \mapright{g_2} j, \quad i \mapright{g_3} j,
\eeq
with $g_1,g_2,g_3 \in G$ distinct, or then arise from a balanced triangle. But given two adjacent edges, there is only one way of completing them into a balanced triangle, while there are at least two ways of completing a 2-subset $i \mapright{g_1} j$, $ i \mapright{g_2} j$ into a 3-subset of $C_n$ of type (\ref{vac}) since $|G| \geq 4$.

This allows to partition $W$ into subsets of this form. To isolate three subsets of the form $W_{ij},W_{ik},W_{jk}$ (with $i,j,k$ distinct), it suffices to take elements $x,y,z$ belonging to different blocks of the partition, such that $xyz \notin H$ (because the only way of this happening is to pick a balanced triangle). Hence, without loss of generality, we may assume that we have identified subsets $W_{12},W_{13},W_{23}$ of $W$ corresponding to elements of $C_n$ of the form $1 \mapright{g} 2$, $1 \mapright{g} 3$ and $2 \mapright{g} 3$, respectively (where $g$ takes all possible values in $G$).

Now we consider an injective mapping $\alpha: W_{12} \rightarrow C_n$ where the image are the elements $x_i$ given by $1 \mapright{g_i} 2$, with $G = \{ g_1,\ldots,g_m\}$. Similarly, we consider an injective mapping $\beta: W_{23} \rightarrow C_n$ where the image are the elements $y_i$ given by $2 \mapright{h_i} 3$, where $G = \{ h_1,\ldots,h_m\}$. But we ignore which elements of $G$ are the $g_i$ or the $h_i$. With these assumptions, for each $i = 1,\ldots,m$, we may identify the element $w_i \in W_{13}$ which should correspond to $1 \longmapright{g_ih_1} 3$: since only 3 parallel edges or balanced triangles are excluded from $\wh{H}_n(G)$, $w_i$ is the unique $w \in W_{13}$ such that $\{ w, \alpha\inv(x_i), \beta\inv(y_1) \} \notin H$. Hence we can define an injective mapping $\gamma:W_{13} \rightarrow C_n$ by associating $\gamma(w_i)$ with $1 \longmapright{g_ih_1} 3$. 

Given $i \in \un{m}$, let $\delta(i) \in \un{m}$ satisfy $g_1h_i = g_{\delta(i)}h_1$. This defines a permutation $\delta$ of $\un{m}$. Moreover, $H$ determines $\delta$ with respect to $\alpha$ and $\beta$ since $w_{\delta(i)}$ is the unique $w \in W_{13}$ such that
$\{ \alpha\inv(g_1), \beta\inv(h_i), w \} \notin H$. Thus $h_i = g_1\inv g_{\delta(i)}h_1$. For all $i,j \in \un{m}$, let $\epsilon(i,j) \in \un{m}$ be such that $g_ih_{\delta\inv(j)} = g_{\epsilon(i,j)}h_1$. Note that the mapping $\epsilon$ is fully determined by $H$ since $w_{\epsilon(i,j)}$ is the unique $w \in W_{13}$ such that $\{ \alpha\inv(g_i), \beta\inv(h_{\delta\inv(j)}), w \} \notin H$. Now we get
$$g_ig_1\inv g_j = g_ig_1\inv g_{\delta(\delta\inv(j))}h_1h_1\inv = g_ih_{\delta\inv(j)}h_1\inv = g_{\epsilon(i,j)},$$
hence $H$ uniquely determines a binary operation $\circ$ on $G = \{ g_1,\ldots,g_m\}$ given by
$$g_i \circ g_j = g_ig_1\inv g_j = g_{\epsilon(i,j)}.$$
But
$$\begin{array}{rcl}
(G,\circ)&\to&(G,\cdot)\\
g&\mapsto&gg_1\inv
\end{array}$$
is obviously an isomorphism, and it is easy to see that this reconstruction of $G$ is essentially unique by symmetry arguments. Therefore
the structure of $G$ can be recovered from $\wh{\H}_n(G)$ as claimed.

(iii) $\Rw$ (i). Immediate.
\qed

\section{Comparing the matroids of Dowling and Rhodes lattices in the general case}

Let $\M_n(G)$ denote the matroid defined by the Dowling lattice $Q_n(G)$. This is a geometric lattice, so $\M_n(G)$ can be defined taking as points the atoms of $Q_n(G)$. These atoms are well known and fall into two types. For all distinct $i,j \in \un{n}$, let $\eta_{ij}$ be the partition of $\un{n}$ with blocks $\{ ij \}$ and $\{ k \}$ (for every $k \in \un{n} \setminus \{ i,j \}$). Write
\bi
\item
$B'_n = \{ (\un{n} \setminus \{ i \}, \iota,[f]_{\iota}) \mid i = 1,\ldots,n;\, f \in F(\un{n} \setminus \{ i \},G) \}$, 
\item
$C'_n = \{ (\un{n}, \eta_{ij} ,[f]_{\eta_{ij}}) \mid 1 \leq i < j \leq n;\,  f \in F(\un{n},G) \}$,
\item
$A'_n = B'_n \cup C'_n$. 
\ei
Then $A'_n$ is the set of atoms of $Q_n(G)$.

Let 
$$\M^{(1)}_n(G) = (B'_n,M^{(1)}_n(G)) = \M_n(G)|_{B'_n} \quad \mbox{ and } \M^{(2)}_n(G) = (C'_n,M^{(2)}_n(G)) = \M_n(G)|_{C'_n}.$$
We intend to compare these matroids with those obtained by the decomposition of the Rhodes matroids (see Lemmas \ref{dsum}(i) and \ref{dsumT}(i)).

By Lemmas \ref{dsum}(ii) and \ref{dsumT}(ii), we have  
$$\H^{(1)}_n(G) \cong U_{n,n} \cong \wh{H}^{(1)}_n(G).$$
It is straightforward to check that also 
$$\M^{(1)}_n(G) \cong U_{n,n}.$$
If $G$ is trivial, we can prove the following:

\bp
\label{twoeq}
For every $n \geq 1$, $\M^{(2)}_n(1) \cong \H^{(2)}_n(1)$.
\ep

\proof
In Section \ref{vlk}, we introduced the operator $\Gamma$ to identify the elements of $C_n$ with group-invertible edges of the form $i \mapright{g} j$, for $i,j \in \un{n}$ distinct and $g \in G$. Since $G$ is trivial, $\Gamma_0$ suffices, hence we view the elements pf $C_n$ as edges of the form $i \edge j$. Then Theorem \ref{thefaces} identifies the faces of $\H^{(2)}_n(1)$ with the forests made with these edges.

On the other hand, since the third component is irrelevant when $G$ is trivial, we can also identify $(\un{n}, \eta_{ij} ,[f]_{\eta_{ij}}) \in C'_n$ with the edge $i \edge j$. Thus we have a natural bijection between $C_n$ and $C'_n$. Moreover, the faces of $\M^{(2)}_n(1)$ are also identified as forests \cite{Zaslavsky}, therefore we obtain an isomorphism between the matroids $\M^{(2)}_n(1) \cong \H^{(2)}_n(1)$.
\qed

Note that Proposition \ref{twoeq} and $\M^{(1)}_n(1) \cong U_{n,n} \cong \H^{(1)}_n(1)$ do not suffice to imply that $\M_n(1) \cong \H_n(1)$
because there is no analogue of Lemma \ref{dsum}(i) in the Dowling world. Indeed, it follows from Theorem \ref{thefaces} that 
$\H_n(1)$ has rank $2n-1$, but $\M_n(1)$ is known to have rank $n$.

From now on, we assume that $n,|G| > 1$. The elements $(\un{n}, \eta_{ij} ,[f]_{\eta_{ij}})\in
C'_n$ can also be efficiently represented through  labeled graphs of the form
$$i \mapright{g} j,$$
where $g = (f(i))\inv f(j) \in G$. Recall that this representation is also used for the atoms of $\wh{\H}^{(2)}_n(G)$, even though they are completely different as SPCs. 

It is well known \cite{Zaslavsky2} that $\M^{(2)}_n(G)$ is the frame matroid of the gain graph $\Delta_n(G)$ defined in Section \ref{alift}. To compare the circuits of the matroids $\M^{(2)}_n(G)$ and $\wh{\H}^{(2)}_n(G)$ (which differ only in one type), the reader is now referred to Section \ref{frali}, and the analysis of frame and lift matroids.

We note also that if $|G| > 1$, then $\M_n(G)$ can be viewed as the frame matroid of the gain graph $\Delta'_n(G)$ obtained by adjoining to each vertex of $\Delta_n(G)$ a loop labeled by some element $g \in G \setminus \{ 1 \}$
\cite{Zaslavsky2}.

\section{Minimal representations}
\label{mini}

Let $\H = (V,H)$ be a BRSC. We remarked in Section \ref{dar} that $\H$ can be represented by $\vee$-generated lattices much smaller than the canonical lattice of flats (but not necessarily geometric, even though if $\H$ is a matroid). More precisely, a lattice representation of $\H$ is a pair $(L,\p)$, where $L$ is a lattice with bottom element $B$, $\p:V \to L\setminus \{ B\}$ a mapping such that $\p(V)$ $\vee$-generates $L$, and $H$ is the set of all subsets of $V$ admitting an enumeration $x_1,\ldots,x_k$ such that
$$\p(x_1) < (\p(x_1) \vee \p(x_2) < \ldots (\p(x_1) \vee \ldots \vee \p(x_k)).$$
This is equivalent to say that $\p(x_1),\ldots,\p(x_k)$ is a transversal of the successive differences for some chain in $L$.
%
%
%
%

We denote by $\LR(\H)$ the class of all
lattice representations of $\H$. Up to isomorphism, every such lattice
representation may be viewed as a sublattice of $\flatx\H$, which
plays then the canonical role of being the largest lattice
representation.

Following \cite[Section 5.4]{brsc},
we may define an ordering on lattice representations of $\H$
by $(L,\p) \geq (L',\p')$ if there exists a $\vee$-map
(i.e. preserving arbitrary joins) $\theta:L \to L'$ such that the diagram
$$\xymatrix{
&A \ar[dl]_{\p} \ar[dr]^{\p'} & \\
L \ar[rr]_{\theta} && L'
}$$
commutes. 
We say that $(L,\p) \in \LR(H)$ is {\em minimal} if 
$$(L,\p) \geq (L',\p') \mbox{ implies $(L,\p) \cong (L',\p')$ for every }(L',\p') \in \LR(H).$$

An element of a lattice $L$ is {\em meet irreducible} if and only if it is covered by a single element.
By \cite[Proposition 5.5.13]{brsc}, $(L,\p) \in \LR(H)$ is minimal if and only if, by identifying a meet irreducible element of $L$ with its unique cover, we never get a lattice representation of $\H$.

We build in this section minimal lattice representations for $\H_n(G)$ and $\M_n(1)$.

Let 
$$L_n = \{ (\emptyset,i) \mid i = 0,\ldots,n-1 \} \cup \{ (I,n) \mid I \subseteq \un{n-1} \}$$
be ordered anti-lexicographically, with inclusion in the first component and the usual order in the second component. It is immediate that $L_n$ is a lattice with
$$((I,i) \vee (J,j)) = (I \cup J, \max\{ i,j\}) \quad \mbox{and} \quad ((I,i) \wedge (J,j)) = (I \cap J, \min\{ i,j\}).$$
We define a mapping $\p:A_n \to L_n$ as follows:
$$\p(\{ i \}, \omega,[f]_{\omega}) = (\emptyset,i) \mbox{ for }i = 1,\ldots,n, \quad 
\p(\{ i,j\}, \omega,[f]_{\omega}) = (\{ i,j\} \setminus \{ n\}, n) \mbox{ whenever } 
1 \leq i < j \leq n.$$
It is immediate that $\p(A_n)$ $\vee$-generates $L_n$. We show that:

\bp
\label{mrh}
$(L_n,\p)$ is a minimal lattice representation of $\H_n(G)$.
\ep

\proof
Assume for the moment that $G$ is trivial (so the third component in $A_n$ is irrelevant and may be omitted).

It follows from Lemma \ref{dsum} and Theorem \ref{thefaces} that the bases of $\H_n(1)$ are of the form $B_n \cup Z$, where $Z  \subseteq C_n$ is such that $\Gamma_0(Z)$ is a spanning tree of $K_n$. Then there exists an ordering 
$n = k_1, k_2, \ldots,k_n$ of $\un{n}$ such that $k_j$ is adjacent in $\Gamma_0(Z)$ to some vertex $k_{i_j} \in \{ k_1,\ldots,k_{j-1} \}$ for $j = 2,\ldots,n$. It is easy to check that $(\{ 1 \}, \omega), \ldots,(\{ n \}, \omega), (\{ k_1,k_2 \}, \omega), (\{ k_{i_3},k_3 \}, \omega), \ldots, (\{ k_{i_n},k_n \}, \omega)$
is a transversal of the successive differences for the chain
\beq
\label{mrh1}
(\emptyset,0) < \ldots < (\emptyset,n) < (\{ k_2\},n) < (\{ k_2,k_3\},n) <\ldots < (\{ k_2,\ldots,k_n\},n)
\eeq
in $L_n$. Thus $(L_n,\p)$ recognizes every basis of $\H_n(1)$ and therefore every independent subset.

Conversely, we consider a transversal of the successive differences for a chain in $L_n$. Since every chain can be refined to a maximal chain, we may assume that our chain is of the form (\ref{mrh1}) for some enumeration $k_2,\ldots,k_n$ of $\un{n-1}$. Write $k_1 = n$. It is easy to check that:
\bi
\item
the unique $x \in C_n$ such that $\p(x) \leq (\emptyset,i)$ and $\p(x) \not\leq (\emptyset,i-1)$ is $(\{ i \}, \omega)$;
\item
the unique $x \in C_n$ such that $\p(x) \leq (\{ k_2\},n)$ and $\p(x) \not\leq (\emptyset,n)$ is $(\{ n,k_2 \}, \omega)$;
\item
the $x \in C_n$ such that $\p(x) \leq (\{ k_2,\ldots,k_i \},n)$ and $\p(x) \not\leq (\{ k_2,\ldots,k_i,n)$ are those of the form $(\{ k_j,k_i \}, \omega)$ with $j < i$.
\ei
This defines precisely a basis of $\H_n(1)$ since it is the union of $B_n$ with the elements corresponding to a spanning tree of $\un{n}$. Therefore $(L_n,\p)$ is a lattice representation of $\H_n(1)$.

Now the meet irreducible elements of $L_n$ are
\beq
\label{meeti}
\{ (\emptyset,i) \mid i = 0,\ldots,n-1 \} \cup \{ (\un{n-1} \setminus \{ i\},n) \mid i \in \un{n-1} \}.
\eeq
If we identify $(\emptyset,i)$ with its unique cover $(\emptyset,i+1)$, we get a lattice of height $2n-2$ which can lo longer be 
a lattice representation of $\H_n(1)$ which has rank $2n-1$.

Out of symmetry, it suffices to show that the lattice $L'$ obtained by identifying $(\un{n-2},n)$ with its unique cover $(\un{n-1},n)$ (and the induced mapping $\p':A_n \to L'$) is no longer a lattice representation of $\H_n(1)$. Let 
$$X = B_n \cup \{ (\{ n,1 \}, \omega), (\{ 1,2 \}, \omega), (\{ 2,3 \}, \omega), \ldots, (\{ n-2,n-1 \}, \omega)\}.$$
Since
$$n \edge 1 \edge 2 \edge 3 \edge \ldots \edge n-1$$
is a spanning tree of $\un{n}$, $X$ is a basis of $\H_n(1)$. Now it is straightforward to check that $X$ is a transversal of the successive differences for a unique chain in $L_n$, namely
$$(\emptyset,0) < \ldots < (\emptyset,n) < (\un{1},n) < (\un{2},n) <\ldots < (\un{n-2},n) < (\un{n-1},n).$$
On the other hand, every chain in $L'$ is also a chain in $L_n$. So
by identifying $(\un{n-2},n)$ with $(\un{n-1},n)$, we lose the only chance we had of recognizing $X$ and therefore $(L',\p')$ is no longer a lattice representation of $\H_n(1)$.

Now it is easy to see that, in view of Lemma \ref{dsum} and Theorem \ref{thefaces}, the same argument can be applied to show that $(L_n,\p)$ is still a minimal lattice representation of $\H_n(G)$ for arbitrary $G$: since we are dealing with forests in $\H_n^{(2)}(G)$, the third components do not play any role because there is never the risk of conflict. 
\qed

Next we build a minimal lattice representation for $\M_n(1)$. Recall than $Q_n(1)$ is isomorphic to the lattice $\Pi_{n+1}$ of full partitions of $\{ 0,\ldots,n\}$, where $\pi \leq \tau$ if and only if every block of $\tau$ is a union of blocks of $\pi$. The set of atoms of this lattice is obviously 
$$E_{n+1} = \{ \eta_{ij} \mid 0 \leq i < j \leq n\},$$
and we may associate each partition $\eta_{ij}$ with the edge $i \edge j$ of $K_{n+1}$ (with vertices $0,\ldots,n$). Then the independent sets correspond to the subforests of $K_{n+1}$ (and the spanning trees are the bases).

Let $2^n$ be the lattice of subsets of $\un{n}$ (ordered by inclusion). We define a mapping $\p':E_{n+1} \to 2^n$ by
$$\p'(\eta_{ij}) = \{ i,j\} \setminus \{ 0\}.$$
We can extend $\p'$ to a mapping $\Phi': \Pi_{n+1} \to 2^n$ by setting $\Phi'(\pi)$ to be the union of the nonsingular blocks of $\pi$ with $0$ removed. It is easy to check that $\Phi'$ is a $\vee$-map and so $\p'(E_{n+1})$ $\vee$-generates $2^n$. A straightforward simplification of the proof of Proposition \ref{mrh} yields

\bp
\label{mrm}
$(2^n,\p')$ is a minimal lattice representation of $\M_n(1)$.
\ep

\section{Boolean matrix representations of minimum degree}
\label{rep}

We have defined BRSC through chains in $\vee$-generated lattices, but an alternative characterization may be provided through Boolean matrices. Given a Boolean matrix $M$ with column space $V$, we define a simplicial complex $\H(M) = (V,H(M))$ as follows. Given $W \subseteq V$, we have $W \in H(M)$ if $M$ admits a (square) submatrix $M[R,W]$ congruent to some lower unitriangular matrix
$$\left(
\begin{matrix}
1&&0&&0&&\ldots&&0\\
?&&1&&0&&\ldots&&0\\
?&&?&&1&&\ldots&&0\\
\vdots&&\vdots&&\vdots&&\ddots&&\vdots\\
?&&?&&?&&\ldots&&1
\end{matrix}
\right)
$$
Two matrices are {\em congruent} if we can transform one into the other
by independently permuting rows/columns. 

By \cite[Corollary 5.2.7]{brsc}, $\H(M)$ is a BRSC. Moreover, every BRSC is of this form, and the matrix $M$ is then said to be a {\em Boolean matrix representation} of $\H$. Given a BRSC $\H$, we denote by $\mindeg(\H)$ the minimum number of rows in a Boolean matrix representation of $\H$.

The following result provides an upper bound for $\mindeg(\H)$.

\bt
\label{mrmd}
Let $(L,\p)$ be a lattice representation of a BRSC $\H$. If $L$ has $p$ meet irreducible elements, then {\rm mindeg}$(\H) \leq p$.
\et

\proof
Assume that $\H = (V,H)$ and $I$ is the set of meet irreducible elements of $L$. Note that $I$ does not contain the top element of $L$ (which is the meet of the empty set). We define an $I \times V$ Boolean matrix $M = (m_{iv})$ by
$$m_{iv} = \left\{
\begin{array}{ll}
0&\mbox{ if }\p(v) \leq i\\
1&\mbox{ otherwise}
\end{array}
\right.
$$
We only need to show that $M$ is a Boolean matrix representation of $\H$.

For each $i \in I$, let $Z_i = \{ v \in V \mid \p(v) \leq i \}$. Suppose that $X \in H \cap 2^{Z_i}$ and $p \in V \setminus Z_i$. Then $\p(X)$ is a transversal of the successive differences for some chain
$a_0 < a_1 < \ldots < a_n$ in $L$. Since $X \subseteq Z_i$, it follows that $\p(X)$ is also a transversal of the successive differences for the chain
$(a_0 \wedge i) < (a_1 \wedge i) < \ldots < (a_n \wedge i)$. Since $\p(p) \not\leq i$, then $\p(X \cup \{ p \})$ a transversal of the successive differences for the chain
$(a_0 \wedge i) < (a_1 \wedge i) < \ldots < (a_n \wedge i) < (a_n \vee \p(p))$ and so $X \cup \{ p\} \in H$. Thus $Z_i \in \flatx\H$.

It follows that $M$ is a submatrix of the $\flatx\H \times V$ Boolean matrix $\mat\H = (q_{Fv})$ defined by
$$q_{Fv} = \left\{
\begin{array}{ll}
0&\mbox{ if }v \in F\\
1&\mbox{ otherwise}
\end{array}
\right.
$$
By \cite[Theorem 5.2.5]{brsc}, $\mat\H$ is a Boolean matrix representation of $\H$, hence $H(\mat(\H)) = H$. Since $M$ is a submatrix of $\mat\H$, we get $H(M) \subseteq H(\mat(\H)) = H$.

Conversely, let $X \in H$. Then there exists an enumeration $x_1,\ldots,x_n$ of the elements of $X$ such that $\p(x_1),\ldots,\p(x_n)$
is a transversal of the successive differences for some chain
$a_0 < a_1 < \ldots < a_n$ in $L$. We may assume that $a_0$ is not the top element (for the case $n = 0$). We claim that $\p(x_1),\ldots,\p(x_n)$
is also a transversal of the successive differences for some chain
\beq
\label{mrmd1}
(i_0 \wedge \ldots \wedge i_n) < (i_1 \wedge \ldots \wedge i_{n}) < \ldots < (i_{n-1} \wedge i_n) < i_n,
\eeq
where $i_0,\ldots,i_n \in I$ and $a_0 \leq (i_0 \wedge \ldots \wedge i_n)$. We use induction on $n$.

The case $n = 0$ holds trivially (the empty set is a transversal of the successive differences for every one element chain and $a_0$ ia a meet of elements from $I$, so $a_0 \leq i_0$ for some $i_0 \in I$).
Assume now that $n \geq 1$ and the claim holds for $n-1$ elements. By the induction hypothesis, $\p(x_2),\ldots,\p(x_n)$
is a transversal of the successive differences for some chain
$$(i_1 \wedge \ldots \wedge i_n) < (i_2 \wedge \ldots \wedge i_{n}) < \ldots < (i_{n-1} \wedge i_n) < i_n,$$
where $i_1,\ldots,i_n \in I$ and $a_1 \leq (i_1 \wedge \ldots \wedge i_n)$. Now $a_0 = (j_1 \wedge \ldots \wedge j_m)$ for some $j_1, \ldots,j_m \in I$. Since $\p(x_1) \not\leq a_0$, it follows that $\p(x_1) \not\leq j_q$ for some $q \in \un{m}$. Let $i_0 = j_q$. Then $\p(x_1) \not\leq (i_0 \wedge \ldots \wedge i_n)$. Since $\p(x_1) \leq a_1 \leq (i_1 \wedge \ldots \wedge i_n)$, we get $(i_0 \wedge \ldots \wedge i_n) < (i_1 \wedge \ldots \wedge i_n)$. Hence $\p(x_1),\ldots,\p(x_n)$
is a transversal of the successive differences for (\ref{mrmd1}), with $i_0,\ldots,i_n \in I$. Finally, $a_0 \leq i_0$ and $a_1 \leq (i_1 \wedge \ldots \wedge i_n)$ together yield $a_0 \leq (i_0 \wedge \ldots \wedge i_n)$, completing the proof of our claim.

Now by permuting the rows and columns of $M$, we may assume that $i_{n-1} < \ldots < i_1 < i_0$ as rows of $M$ and $x_n < \ldots < x_2 < x_1$ as columns of $M$. We claim that the submatrix $M[i_{n-1},\ldots,i_0,x_n,\ldots,x_1]$ is lower unitriangular. This amounts to check that
$m_{i_{k},x_{j}} = 0$ whenever $1 \leq j \leq k \leq n-1$
and $m_{i_{j-1},x_{j}} = 1$ for $j \in \un{n}$. Indeed, we have $\p(x_j) \leq (i_j \wedge \ldots \wedge i_{n}) \leq i_k$ and $\p(x_j) \not\leq (i_{j-1} \wedge \ldots \wedge i_{n})$, in view of $\p(x_j) \leq (i_j \wedge \ldots \wedge i_{n})$, yields $\p(x_j) \not\leq i_{j-1}$. Hence $M[i_{n-1},\ldots,i_0,x_n,\ldots,x_1]$ is lower unitriangular and so $X \in H(M)$. Thus $H \subseteq H(M)$ and so $H = H(M)$. Therefore $M$ is a Boolean matrix representation of $\H$ and the proof is complete.
\qed

We can now apply this result to the matroids $\H_n(1)$ and $\M_n(1)$.

\bc
\label{mindeg}
${\rm mindeg}(\H_n(1)) = 2n-1$ and ${\rm mindeg}(\M_n(1)) = n$.
\ec

\proof
The meet irreducible elements of $L_n$ were computed in (\ref{meeti}), and there are $2n-1$ of them. Now Proposition \ref{mrh} and Theorem \ref{mrmd} yield ${\rm mindeg}(\H_n(1)) \leq 2n-1$. On the other hand, $\H_n(1)$ has rank $2n-1$, so we need at least $2n-1$ rows in a matrix representation to recognize the bases. Therefore ${\rm mindeg}(\H_n(1)) = 2n-1$.

Similarly, ${\rm mindeg}(\M_n(1)) = n$ since $\M_n(1)$ has rank $n$ and the lattice $2^n$ (recall Section \ref{mini}) has $n$ meet irreducible elements (of the form $\un{n} \setminus \{ i\}$, for $i \in \un{n}$).
\qed

We should note that the minimum degree being equal to the rank is by no means standard behavior. Even if $\H$ is a matroid of rank 2. Take for instance the uniform matroid $U_{2,2^n-1}$ (that is, the complete graph $K_{2^n-1}$ viewed as a simplicial complex). Then $U_{2,2^n-1}$ is recognized by a $n \times (2^n-1)$ Boolean matrix having all columns distinct and nonzero. However, any matrix with less than $n$ rows would have necessarily two equal columns and would fail to recognize a 2-set. Therefore $\mindeg(U_{2,2^n-1}) = n$, even though $U_{2,2^n-1}$ has rank 2.

\section{Open problems}

We computed in Sections \ref{mini} and \ref{rep} minimal latttice representations and Boolean representations of minimum degree for the Dowling and Rhodes matroids, in the case of the group being trivial. Can this be done for an arbitrary finite group?

\smallskip

Can this study be generalized to the frame and lift matroids of arbitrary gain graphs, or even to the frame and lift matroids of arbitrary biased graphs?

\section*{Acknowledgments}

The first author acknowledges support from the Binational Science Foundation (BSF) of the United States and Israel, grant number 2012080. The second author acknowledges support from the Simons Foundation (Simons Travel Grant Number 313548).
The third author was partially supported by CMUP (UID/MAT/00144/2013), which is funded by FCT (Portugal) with national (MEC) and European structural funds (FEDER), under the partnership agreement PT2020.


\bigskip

{\sc Stuart Margolis, Department of Mathematics, Bar Ilan University,
  52900 Ramat Gan, Israel}

{\em E-mail address:} margolis@math.biu.ac.il

\bigskip

{\sc John Rhodes, Department of Mathematics, University of California,
  Berkeley, California 94720, U.S.A.}

{\em E-mail addresses}: rhodes@math.berkeley.edu, BlvdBastille@gmail.com

\bigskip

{\sc Pedro V. Silva, Centro de
Matem\'{a}tica, Faculdade de Ci\^{e}ncias, Universidade do
Porto, R. Campo Alegre 687, 4169-007 Porto, Portugal}

{\em E-mail address}: pvsilva@fc.up.pt

\end{document}